\documentclass[twoside,11pt]{article}
\title{} \author{} \date{}
\usepackage{latexsym,amssymb,times}
\input amssym.def
\newtheorem{te}{Theorem}[section]
\newtheorem{prop}[te]{Proposition}

\newtheorem{fac}[te]{Fact}

\newtheorem{lem}[te]{Lemma}

\newtheorem{ex}[te]{Example}

\def\dok{\noindent{\bf Proof. }}
\def\kdok{\hfill $\Box$ \par \vspace*{2mm} }
\def\a{\alpha}

\def\f{\varphi}
\def\p{\psi}
\def\o{\omega}
\def\k{\kappa}
\def\l{\lambda}

\def\r{\rho}
\def\s{\sigma}

\def\t{\tau}

\def\T{{\mathbb T}}

\def\Q{{\mathbb Q}}
\def\B{{\mathbb B}}
\def\N{{\mathbb N}}
\def\X{{\mathbb X}}
\def\Y{{\mathbb Y}}
\def\Z{{\mathbb Z}}

\def\A{{\mathbb A}}
\def\C{{\mathbb C}}

\def\BH{{\mathbb H}}
\def\BK{{\mathbb K}}
\def\BI{{\mathbb I}}

\def\BR{{\mathbb R}}

\def\CC{{\mathcal C}}

\def\CT{{\mathcal T}}

\def\K{{\mathcal K}}

\def\CF{{\mathcal F}}

\def\c{{\mathfrak{c}}}

\def\la{\langle}
\def\ra{\rangle}

\def\dom{\mathop{\mathrm{dom}}\nolimits}
\def\ran{\mathop{\mbox{ran}}\nolimits}
\def\id{\mathop{\mathrm{id}}\nolimits}
\def\rank{\mathop{\rm rank}\nolimits}
\def\Card{\mathop{\rm Card}\nolimits}
\def\Iso{\mathop{\rm Iso}\nolimits}
\def\Aut{\mathop{\rm Aut}\nolimits}
\def\Sym{\mathop{\rm Sym}\nolimits}

\def\ar{\mathop{\rm ar}\nolimits}
\def\Mod{\mathop{\rm Mod}\nolimits}

\def\Var{\mathop{\rm Var}\nolimits}

\def\LO{\mathop{\rm LO}\nolimits}

\def\tp{\mathop{\rm tp}\nolimits}

\def\Th{\mathop{\rm Th}\nolimits}

\def\Pa{\mathop{\rm Pa}\nolimits}
\def\pred{\mathop{\rm pred}\nolimits}

\def\Lev{\mathop{\rm Lev}\nolimits}

\def\bcd{\dot{\bigcup}}

\def\du{\mathrel{\dot{\cup}}}
\begin{document}
\thispagestyle{plain}
\begin{center}
           {\large \bf \uppercase{Sharp Vaught's Conjecture for some classes of partial orders}}
\end{center}
\begin{center}
{\bf Milo\v s S.\ Kurili\'c}\footnote{Department of Mathematics and Informatics, Faculty of Sciences, University of Novi Sad,
                                      Trg Dositeja Obradovi\'ca 4, 21000 Novi Sad, Serbia.
                                      e-mail: milos@dmi.uns.ac.rs}
\end{center}
\begin{abstract}
\noindent
Matatyahu Rubin has shown that a sharp version of Vaught's conjecture, $I(\CT,\o)\in \{ 0,1,\c\}$, holds for each complete theory of linear order $\CT$.
We show that the same is true for each complete theory of partial order
having a model in the minimal class of partial orders containing the class of linear orders
and which is closed under finite products and finite disjoint unions.
The same holds for the extension of the class of rooted trees admitting a finite monomorphic decomposition, obtained in the same way.
The sharp version of Vaught's conjecture also holds for the theories of trees which are infinite disjoint unions of linear orders.


\end{abstract}
\section{Introduction}\label{S1}
We recall that Vaught's conjecture \cite{Vau}, stated by Robert Vaught in 1959,
is the statement that the number $I(\CT ,\o)$ of non-isomorphic countable models
of a complete countable first-order theory $\CT$
is either at most countable or continuum.
The long list of the results related to this (still open) problem includes its reductions:
Arnold Miller has shown that the full conjecture follows from its restriction to complete theories of partial orders, see \cite{Stee},
while the reduction to lattices was proved by Taitslin \cite{Tai} (see also \cite{Hodg}, p.\ 231).

Vaught's conjecture was confirmed for several classes of partial orders:
for linear orders with countably many unary predicates by Rubin \cite{Rub},
for (model-theoretic) trees by Steel \cite{Stee},
for arborescent structures, in particular, for the class of partial orders not embedding the four element poset $N$
(containing trees and reversed trees) by Schmerl \cite{Sch3},
for Boolean algebras by Iverson \cite{Ive}.

Generally speaking, if Vaught's conjecture is true for a theory $\CT$,
it is natural to ask what is going on
with the theories of reducts of definitional expansions of models of $\CT$.
Following that concept, on the basis of the results of
Rubin \cite{Rub}, concerning colored linear orders;
Fra\"{\i}ss\'{e}, about monomorphic and almost monomorphic structures (see \cite{Fra});
and Gibson, Pouzet and Woodrow \cite{Gib} describing (roughly speaking) indiscernible sequences for a structure;
Vaught's conjecture was confirmed for monomorphic theories in \cite{KMon} and for almost chainable theories in \cite{KAC}.
In \cite{KFMD}, using the results of Pouzet and Thi\'{e}ry \cite{PT} concerning the structures admitting a finite monomorphic decomposition (FMD structures),
Vaught's conjecture was confirmed for FMD theories (theories of reducts of quantifier-free definitional expansions of linear orders colored in finitely many convex colors).

In fact, Rubin has shown that a sharp version of Vaught's conjecture, VC$^\sharp$: $I(\Th (\X),\o)\in \{ 0,1,\c\}$,
is true whenever $\X$ belongs to the class LO of linear orders.
In addition by \cite{Ive} and \cite{KFMD}, VC$^\sharp$ is true for complete theories of Boolean algebras, for FMD-theories,
and, in particular, for theories of FMD posets.

In this paper we consider partial orders and the corresponding relational language $L_b=\la R\ra$, where $\ar (R)=2$.
One more idea to extend a class of structures satisfying VC is to regard its closure under some operations.
So, after preliminaries given in Section \ref{S2},
in Section \ref{S3} for an $\cong$-closed class $\K$ of $L_b$-structures we describe its minimal closure under finite products and finite disjoint unions
and show that, under certain additional conditions, the VC$^\sharp$ for (the theories of) the structures from $\K$ extends to the VC$^\sharp$ for the structures from $\la \K\ra$.
In Sections \ref{S4}--\ref{S6} we confirm  VC$^\sharp$ for the theories of the partial orders from  $\la \K \ra$,
where $\K$ is the class of (a) initially finite trees satisfying VC$^\sharp$; (b) rooted FMD trees; (c) linear orders.
In Section \ref{S7} we show that VC$^\sharp$ holds for the  trees which are infinite disjoint unions of linear orders.
(This class is exactly the intersection of the class of trees and the class of reversed trees, for which Vaught's conjecture was already confirmed by Steel.)

For convenience, for a complete theory $\CT$ we write $I(\CT,\o)= 0$ iff $\CT$ has a finite model;
we will say that $\CT$ satisfies VC$^\sharp$ iff $I(\CT,\o)\in\{0, 1,\c\} $.
Also,  $\CT$ is called $\o$-categorical iff $I(\CT,\o)\in\{ 0,1\} $.
This terminology is extended to structures in the usual way: e.g.\ $\X$ satisfies VC$^\sharp$ iff $\Th(\X)$ does.
For notational simplicity instead of $I(\CT, \o)$ and $I(\Th (\X), \o)$ we will write only $I(\CT)$ and $I(\X )$.
For an $L_b$-structure $\X$ by $\tp (\X )$ we denote the {\it isomorphism type of $\X$} (the class of all $L_b$-structures isomorphic to $\X$).
\section{Products and lexicographical sums of $L_b$-structures}\label{S2}
\paragraph{Direct products}
We will use the following facts, true for any language $L$.
\begin{fac}\label{T043}
Let $\X$, $\Y$, $\Z$, $\X _i$ and $\Y _i$, for $i\in I$, be $L$-structures. Then

(a) $\prod _{i\in I}\X _i \cong \prod _{i\in I}\X _{\pi (i)}$, for any permutation $\pi \in \Sym (I)$;

(b) $\prod _{i\in I}\X _i \cong \prod _{j\in J}\prod _{i\in I _j}\X _i$, for any partition $\{ I_j :j\in J\}$ of $I$;

(c) If $\X _i \cong \Y _i$, for all $i\in I$, then $\prod _{i\in I}\X _i \cong \prod _{i\in I}\Y _i$;

(d) If $\X _i \equiv \Y _i$, for all $i\in I$, then $\prod _{i\in I}\X _i \equiv \prod _{i\in I}\Y _i$;

(e) If $|L|\leq \o$, $|I|<\o$ and $I(\X _i )\leq 1$, for all $i\in I$, then $I(\prod _{i\in I}\X _i)\leq 1$.
\end{fac}
\dok
Statements (a)--(c) are elementary (see \cite{CK}, p.\ 326),
(d) follows from the Feferman-Vaught analysis of reduced products \cite{FV} (see also \cite{CK}, p.\ 345),
while (e) was proved by Grzegorczyk \cite{Grz} (see also \cite{Hodg}, p.\ 347).
\hfill $\Box$
\begin{fac}\label{T045}
If $\X$ and $\Y$ are countable $L$-structures
and $\X \times \Y ' \not\cong \X \times \Y ''$, whenever $\Y '$ and $\Y ''$ are non-isomorphic countable models of $\Th (\Y )$,
then $I(\X  \times \Y ) \geq I(\Y )$.
\end{fac}
\dok
Let $I(\Y )=\k$ and let $\Y _\a$, for $\a <\k$, be non-isomorphic countable models of $\Th (\Y )$.
By Fact \ref{T043}(d) $\X \times\Y _\a$, $\a <\k$, are countable models of $\Th (\X  \times \Y )$
and by our assumption they are non-isomorphic.
\hfill $\Box$
\paragraph{Lexicographical sums of $L_b$-structures}
Let $\BI=\la  I, \r _I\ra$ be an $L_b$-structure and let $\X _i =\la X_i , \r _i \ra$, $i\in I$, be $L_b$-structures with pairwise disjoint domains.
The {\it lexicographical sum of the structures $\X _i$, $i\in I$, over the structure $\BI$}, in notation $\sum _{\BI}\X _i$, is the  $L_b$-structure
$\X :=\la X, \r \ra$, where $X := \bigcup _{i\in I}X_i$ and $\r :=\bigcup _{i\in I}\r _i \cup \bigcup _{\la i,j \ra\in \r _I \setminus \Delta _I} X_i \times X_j$.
Namely, for $x,x'\in X$ we have: $x\,\r \,x'$ iff
\begin{equation}\label{EQ200}
\exists i\in I \;\; (x,x'\in X_i \land x \,\r _i\, x') \;\lor\;
\exists \la i,j \ra\in \r _I \setminus \Delta _I \;\; (x\in X_i \land x'\in X_j ).
\end{equation}
If, in particular,  $\r _I=\emptyset$, then the corresponding lexicographical sum is the structure $\bcd _{i\in I} \X _i :=\la \bigcup _{i\in I} X_i , \bigcup _{i\in I} \r _i\ra$
called the {\it disjoint union} of the structures $\X _i$, $i\in I$.

It is evident that the definition of $\sum _{\BI}\X _i$ does not depend on $\r _I \cap \Delta _I$; so w.l.o.g.\ we can assume that the relation $\r _I$ is
reflexive or irreflexive. It is easy to see that a lexicographical sum of partial orders over a partial order is a partial order
and the same holds for  strict partial orders, tournaments, graphs, etc.
\begin{prop}\label{T216}
Let $\X$ be a countable $L_b$-structure, $\{ X_i : i\in I\}$ a partition of $X$ and
$\X_i =\X \upharpoonright X_i$, $i\in I$, the corresponding substructures of $\X$.
Then we have

(a) If $|I|<\o$ and  $\bigcup _{i\in I}f_i \in \Aut (\X )$, whenever $f_i \in \Aut (\X_i )$, for $i\in I$, then
\begin{equation} \label{EQ237}\textstyle
\forall i\in I \;(\X_i \mbox{ is $\o$-categorical}) \Rightarrow \X  \mbox{ is $\o$-categorical;}
\end{equation}

(b) If for each $F\in \Aut (\X  )$ and $i\in I$ from $F[X_i]\cap X_i \neq \emptyset$ it follows that $F[X_i]=X_i$,
then the implication ``$\Leftarrow$" holds in (\ref{EQ237}).
\end{prop}
\dok
We recall one incarnation of the theorem of Engeler, Ryll-Nardzewski and Svenonius (see \cite{Hodg}, p.\ 341):
{\it a countable structure $\X$ of a countable language is $\o$-categorical iff $|X^n /\sim _n^\X|<\o$, for all $n\in \N$,}
where  $\sim _n^\X$ is  the equivalence relation on the set $X^n$
defined by $\bar x \sim _n^\X \bar y $ iff  $f\bar x =\bar y$, for some $f\in \Aut (\X )$.

(a) Under the assumptions, for each $i\in I$ and $n\in \N$ there are $m^i_n\in \N$ and $\bar x ^{i,j}=\la x ^{i,j}_k:k<n\ra \in X_i^n$, for $j<m^i_{n}$,  such that we have
$X_i ^n /\!\sim _{n}^{\X _i} \; =\{ [\bar x ^{i,j}]_{\sim _{n}^{\X _i} }: j<m^i_{n}\}$.
Let $n\in \N$ and $\bar z=\la z_r :r<n\ra \in X^n$.
Then $n=\bigcup _{i\in I}S_i $, where $S_i:=\{ r<n :z_r \in X_i\}$.

If $n_i:=|S_i| >0$
and $S_i=\{ r^i_k :k <n_i\}$ is the enumeration such that $r^i_0 <\dots <r^i_{n_i -1}$,
then $\la z _{r^i_k}:k<n_i\ra \in X_i ^{n_i}$
and, hence, there is $j_i<m^i_{n_i}$ such that
\begin{equation} \label{EQ218}
\bar z| S_i :=\la z _{r^i_k}:k<n_i\ra \sim _{n_i}^{\X _i}\la x ^{i,j_i}_k:k<n_i\ra =\bar x ^{i,j_i},
\end{equation}
which means that there is $f_i\in \Aut (\X_i)$ such that
\begin{equation} \label{EQ212}
\forall k <n_i \;\;f_i(z _{r^i_k})=x ^{i,j_i}_k.
\end{equation}
Defining $J:=\{ i\in I : n_i >0\}$ we have $n=\bigcup _{i\in J}S_i =\bigcup _{i\in J}\{ r^i_k :k <n_i\}$
and we define the $n$-tuple $\bar y=\la y_r :r<n\ra\in X^n$ by
\begin{equation}\label{EQ215}
y_{r^i_k}  :=  x ^{i,j_i}_k, \mbox{ for } i\in J \mbox{ and } k<n_i.
\end{equation}
For $i\in I\setminus J$ we take $f_i =\id _{X_i}$, by the assumption we have $F:=\bigcup _{i\in I}f_i  \in \Aut (\X )$,
by (\ref{EQ212}) and (\ref{EQ215}) we have $F \bar z=\bar y$
and, hence, $\bar z \sim _n ^\X \bar y$.

Thus to each $\bar z\in X^n$ we can adjoin
the sets $S_i^{\bar z}$, $i\in I$, defined as above;
$j_i^{\bar z} <m^i_{n_i}$, for $i\in J^{\bar z}:=\{ i\in I : S_i^{\bar z}\neq \emptyset \}$, such that (\ref{EQ218}) holds;
and $\bar y ^{\bar z}\sim _n ^\X \bar z$ defined by (\ref{EQ215}).
If for $\bar z,\bar z'\in X^n$ all the aforementioned parameters are the same,
then  we have $\bar y ^{\bar z}=\bar y ^{\bar z'}$ and, hence, $\bar z \sim _n ^{\X }  \bar z'$.
So, since there are finitely many such sequences $\la S_i :i\in I\ra$,
and the sets $X_i /\!\sim _{n}^{\X _i}$ are finite,
there are finitely many codes $\bar y ^{\bar z}$ for the orbits of $n$-tuples $\bar z \in X^n$
and, hence, $X ^n /\!\sim _{n}^{\X } <\o$.
Thus the structure $\X $ is $\o$-categorical.

(b) Let the assumption of (b) hold and let the structure $\X $ be $\o$-categorical.
Suppose that there are $i\in I$, $n\in \N$, and $\bar z^k \in X_i ^n$, for $k\in \o$, such that $\bar z^k \not\sim ^{\X _i}_n\bar z^l$, for $k\neq l$.
Since $\X $ is $\o$-categorical there are $S\in [\o]^\o$ and $\bar y\in X^n$
such that $\bar z^k \sim ^{\X }_n\bar y$, for all $k\in S$.
So, if $k,l\in S$ and $k\neq l$, then $\bar z^k \sim ^{\X }_n\bar z^l$
and, hence, there is $F\in \Aut (\X )$ such that $F\bar z^k =\bar z^l $.
By the assumption and since $\bar z^k,\bar z^l \in X_i^n$ we have  $f:=F|X_i \in \Aut (\X _i)$.
Since $f\bar z^k =\bar z^l$ we have $\bar z^k \sim ^{\X _i}_n\bar z^l$ and we obtain a contradiction.
\hfill $\Box$
\begin{prop}\label{T200}
Let $\sum _{\BI}\X _i$ and $\sum _{\BI}\Y _i$ be lexicographical sums of $L_b$-structures. Then we have

(a) If $\X _i \cong  \Y _i$, for all $i\in I$, then $\sum _{\BI}\X _i \cong  \sum _{\BI}\Y _i $;

(b) If $\X _i \equiv  \Y _i$, for all $i\in I$, then $\sum _{\BI}\X _i \equiv  \sum _{\BI}\Y _i $;

(c) If $|I|<\o$ and $I(\X _i )\leq 1$, for all $i\in I$, then $I(\sum _{\BI}\X _i)\leq 1$.
\end{prop}
\dok
(a) Let $\X:=\sum _{\BI}\X _i =\la X,\r\ra$ and $\Y:=\sum _{\BI}\Y _i =\la Y,\s\ra$.
If $f_i\in \Iso (\X _i,\Y _i)$, for $i\in I$, then, by (\ref{EQ200}), $\bigcup _{i\in I}f_i \in \Iso (\X ,\Y )$.

(b) Clearly, $\{ X_i: i\in I\}$ and $\{ Y_i: i\in I\}$ are partitions of $X$ and $Y$ respectively.
By Remark 1.5.4 of \cite{Stee}, if for different $i,j\in I$, $x\in X_i$, $y\in Y_i$, $x' \in X_j$ and $y '\in Y_j$,
$(\X _i,x)\equiv _0 (\Y _i,y)$ and $(\X _j,x')\equiv _0 (\Y _j,y')$ implies that $ x \mathrel{\r } x' \Leftrightarrow y \mathrel{\s } y' $,
then  for each $\a <\o _1$
\begin{equation}\label{EQ236}
(\forall i\in I \;\;\X _i \equiv _\a \Y _i ) \Rightarrow \X \equiv _\a \Y .
\end{equation}
But if $i$, $j$, $x$, $y$, $x'$ and $y '$ are as assumed,
then by (\ref{EQ200}) we have $x\,\r \,x'$ iff $i \,\r _I\, j$ iff $y\,\s \,y'$. So the implication is true and we have (\ref{EQ236}).

Let $\X _i \equiv  \Y _i$, for all $i\in I$. Then for each $n\in \o$ we have $\X _i \equiv_n  \Y _i$, for all $i\in I$,
and, by (\ref{EQ236}), $\X  \equiv_n  \Y $. Thus $\X \equiv \Y $.

(c) follows from the downward L\"{o}wenheim-Skolem theorem, the proof of (a) and Proposition \ref{T216}(a).
\hfill $\Box$
\paragraph{Disjoint unions of connected $L_b$-structures}
If $\X =\la X,\r \ra$ is an $L_b$-structure, then
the transitive closure $\r _{rst}$ of the relation $\r _{rs} =\Delta _X \cup \r \cup \r ^{-1}$
(given by $x \;\r_{rst} \;y$ iff there are $n\in \N$ and $z_0 =x , z_1, \dots ,z_n =y$ such that $z_i \;\r _{rs} \;z_{i+1}$, for each $i<n$)
is the minimal equivalence relation on $X$ containing $\r$.
The equivalence classes $[x]_{\r _{rst}}$, $x\in X$, are the {\it connected components} of $\X$ and $\X$ is said to be {\it connected}
iff $|X/\r _{rst}|=1$.
Clearly, the formula $\f _{rs}(u,v):= u=v \lor R(u,v) \lor R(v,u)$ defines the relation $\r_{rs}$ in the structure $\X$
and for $n\in \N$ the formula
$\p _n(u,v):=\exists w_0, \dots w_n \;( u=w_0 \land v=w_n \land \bigwedge _{i<n} \f _{rs}(w_i,w_{i+1}))$
says that there is a path of length $\leq n$ from $u$  to $v$.
So,  the relation $\r_{rst}$ is defined in $\X$ by the $L_{\o _1\o}$-formula $\f _{rst} (u,v):=\bigvee _{n\in \N} \p _n(u,v)$
and the $L_{\o _1\o}$-sentence $\f _{conn}:=\forall u,v \;\f _{rst} (u,v)$ says that an $L_b$-structure is connected.
\begin{fac}\label{T226}
For each uncountable connected $L_b$-structure $\X$ there is a countable connected $L_b$-structure $\Y \equiv \X$.
\end{fac}
\dok
For the $L_{\o _1\o}$-sentence $\f :=\bigwedge \Th (\X) \land \f _{conn}$
we have $\X \models \f $ and by the weak downward L\"{o}wenheim-Skolem theorem for $L_{\o _1\o}$ (see \cite{Hodg}, p.\ 91)
there is a countable $L_b$-structure $\Y \models \f$. So, $\Y$ is connected and $\Y \equiv \X$.
\hfill $\Box$
\begin{fac}\label{T4015}
If $\X _i$, $i\in I$, and $\Y _j$, $j\in J$, are families of pairwise disjoint connected $L_b$-structures
and $F:\bcd _{i\in I} \X _i  \rightarrow \bcd _{j\in J} \Y _j $ is an isomorphism,
then there is a bijection $f:I\rightarrow J$ such that $F[X_i]=Y_{f(i)}$ and, hence, $\X_i\cong\Y_{f(i)}$, for all $i\in I$.
\end{fac}
\dok
Let $\X=\bcd _{i\in I} \X _i =\la X,\r\ra$ and $\Y=\bcd _{j\in J} \Y _j =\la Y,\s\ra$.
Then $\r _{rst}$ and $\s _{rst}$ are the equivalence relations
corresponding to the partitions $\{ X _i : i\in I\}$ of $X$ and $\{ Y _j : j\in J\}$ of $Y$.
Since isomorphisms preserve $L_{\o _1\o}$-formulas in both directions
$F:\la X,\r _{rst}\ra \rightarrow \la Y,\s_{rst}\ra$ is an isomorphism.
So, for each $i\in I$ there is (a unique) $j_i\in J$ such that $F[X_i]=Y_{j_i}$ and $f:I\rightarrow J$ given by $f(i)=j_i$ is well defined.
Since $F$ is an injection (surjection) $f$ is an injection (surjection).
\kdok
For brevity, by $(\ast)$ we denote the assumption that $\X _i$, $i\in I$, are pairwise disjoint connected $L_b$-structures and $|I|\leq \o$.
\begin{prop}\label{T240}
Let $(\ast)$ hold and $|\X _i|\leq \o$, for $i\in I$. Then $I(\bcd _{i\in I}\X _i) \leq 1$ iff  $I(\X _i )\leq 1$, for all $i\in I$, and $|\{ \tp (\X _i) :i\in I \}| <\o$.
\end{prop}
\dok
($\Rightarrow$) Let the structure $\X := \dot{\bigcup} _{i\in I}\X _i$ be $\o$-categorical.
By Fact \ref{T4015} if $F\in \Aut (\X  )$ and $F[X_i]\cap X_i \neq \emptyset$, then $F[X_i] = X _i$,
and by Proposition \ref{T216}(b) the structures $\X _i$, for all $i\in I$, are $\o$-categorical.
Suppose that $\{ \tp (\X _i ) :i\in I \} =\{ \t _n:n<\o\}$.
For each $n\in \o$ we take $i_n \in I$ such that $\tp(\X _{i_n})=\t _n$ and pick an $x_n\in \X _{i_n}$;
clearly, $m\neq n$ implies that $i_m \neq i_n$ and $x_m \neq x_n$.
Since $|X /\! \sim ^{\X}_1|<\o$ there are different $m,n\in\o$ such that $x_m \sim ^\X_1 x_n$;
so, there is $F\in \Aut (\X )$ such that $F(x_m)=x_n$.
By Fact \ref{T4015} $\X _{i_m}\cong\X _{i_n}$, that is, $\t_m = \t _n$, which is false.

($\Leftarrow$) Let the structures $\X _i$, $i\in I$, be $\o$-categorical.
If $|I|<\o$, then $\X$ is $\o$-categorical by Proposition \ref{T200}(c).
Otherwise, let $I=\o$ and first suppose that $\tp (\X _i )= \t $, for all $i\in \o$.
Let $n\in \N$ and $\Y :=\bcd _{i<n}\X _i$.
If $\bar x =\la x_0, \dots ,x_{n-1}\ra \in X ^n$,
then there is $K \in [\o ]^n$ such that $\bar x\in (\bigcup _{i\in K}X_i )^n$
and, clearly, there is $F\in \Aut (\X )$ such that $F[\bigcup _{i\in K}X_i]=Y$.
So, $\bar x \sim ^\X_n F \bar x \in Y^n$ and, thus,
\begin{equation}\label{EQ275}\textstyle
\forall \bar x \in X^n \;\;\exists \bar y \in Y^n \;\; \bar x \sim ^\X_n\bar y .
\end{equation}
Suppose that there are $\bar x^k \in X^n$, $k\in \o$, such that $\bar x^k \not\sim ^\X_n \bar x^l$, for $k\neq l$.
By (\ref{EQ275}) there are $\bar y^k \in Y^n$ such that $\bar y ^k \sim ^\X_n \bar x^k$, for $k\in \o$.
By  Proposition \ref{T200}(c)  $\Y $ is $\o$-categorical
and, hence, there are different $k,l\in \o$ such that $\bar y ^k \sim ^\Y_n \bar y^l$;
thus, there is $g\in \Aut (\Y )$ such that $g \bar y ^k=\bar y^l$.
Clearly $h :=g\cup \id _{X\setminus Y} \in \Aut (\X )$.
So, since $h\bar y ^k=g \bar y ^k=\bar y^l$ we have $\bar y^k \sim ^\X_n \bar y^l$
and, hence, $\bar x^k \sim ^\X_n \bar x^l$, which gives a contradiction.

Let $|\{ \tp (\X _i) :i\in I \}|=\{\t _m :m<n  \}$
and let $I_m :=\{ i\in I : \tp (\X_i )= \t _m\}$ and $\X ^m :=\bcd _{i\in I_m}\X _i$, for $m<n$.
Then the structures $\X _m$, $m<n$, are $\o$-categorical
and $\X =\bcd _{m<n}\X _m$ is $\o$-categorical by  Proposition \ref{T200}(c).
\hfill $\Box$

\begin{lem}\label{T042}
If $I$, $J$ and $K$ are pairwise disjoint non-empty sets, where $|I|<\o$,
$f:I\cup J \rightarrow I\cup K$ is a bijection
and $\r$ is an equivalence relation on the set $I\cup J\cup K$ such that $f\subset \r$,
then there is a bijection $\pi : J\rightarrow K$ such that $\pi \subset \r$.
\end{lem}
\dok
Clearly, the set $f =\{ \la s, f(s)\ra :s\in I\cup J\}$ is a binary relation on the set $I\cup J\cup K$.
Let $\,\sim\;:=\,f_{rst}$ be the minimal equivalence relation on the set $I\cup J\cup K$ containing $f$.
We construct a bijection $\pi : J\rightarrow K$ such that
\begin{equation}\label{EQ027}
j \sim \pi (j), \mbox{ for all }j\in J.
\end{equation}
For $j\in J$, let $[j]$ denote the $\sim$-equivalence class of $j$.
First we prove that
\begin{equation}\label{EQ007}
\forall j\in J \;\; f[[j] \cap (I\cup J)]= [j] \cap (I\cup K) .
\end{equation}
If $s\in [j] \cap (I\cup J)$,
then, since $f \subset \;\sim$, we have $f(s)\sim s \sim j$
and, hence, $f(s)\in [j]$.
In addition, $s\in I\cup J=\dom f$,
which implies that $f(s)\in \ran f =I\cup K$.

Conversely, if $s\in [j] \cap (I\cup K)$,
then, since $f$ is a surjection, there is $r\in I\cup J$ such that $f(r)=s$.
Since  $s\in [j]$ we have $r\sim f(r)= s \sim j$,
which gives $r\in [j]$.
Thus $r\in [j] \cap (I\cup J)$ and $s=f(r)\in f[[j] \cap (I\cup J)]$;
so, (\ref{EQ007}) is true.

Since the mapping $f$ is a bijection, $I\cap J=\emptyset$ and $I\cap K=\emptyset$, by (\ref{EQ007})
we have
\begin{equation}\label{EQ008}
f[[j] \cap (I\cup J)]=f[[j] \cap I]\du f[[j] \cap J]=([j] \cap I)\du([j] \cap K)
\end{equation}
and $|I|<\o$ gives $m:=|f[[j] \cap I]|=|[j] \cap I|\in \o$.
Thus, by  (\ref{EQ008}) and since $|f[[j] \cap J]|=|[j] \cap J|$, we have
$m + |[j] \cap J|=m+|[j] \cap K|$.
Now we have:
if $|[j] \cap J|=\k \geq \o$, then $|[j] \cap K|=\k$;
and if $|[j] \cap J|=n\in \N$,  then $|[j] \cap K|=n$.
Thus we have
\begin{equation}\label{EQ010}
\forall j\in J \;\; |[j] \cap J|=|[j] \cap K|.
\end{equation}
If $j\in J$,
then $j\in [j] \cap J$
and by (\ref{EQ010}) we have $|[j] \cap K|=|[j] \cap J|\geq 1$.
Thus for each $j\in J$ there exists $k\in K$ such that $j\sim k$.
We prove that
\begin{equation}\label{EQ005}
\forall k\in K \;\;\exists j\in J \;\; k\sim j .
\end{equation}
Clearly $f^{-1}:I\cup K \rightarrow I\cup J$ is a bijection
and $\sim ':= (f^{-1})_{rst}$ is the transitive closure of the relation $(f^{-1})_{rs}=\Delta \cup f^{-1}\cup f =f_{rs}$.
Thus $\sim ' =\sim$
and, as above, for each $k\in K$ there is $j\in J$ such that $k\sim 'j$,
that is $k\sim j$. So (\ref{EQ005}) is true.

Choosing representatives of $\sim$-classes which intersect the set $J$
we obtain a partition of the set $J$,  $J/\!\sim =\{ [j _t]\cap J :t\in T\}$,
where $j_t\in J$, for all $t\in T$.
By (\ref{EQ005}), for each $k\in K$ there is $t\in T$ such that $k\sim j_t$,
which implies that $K/\sim =\{ [j _t]\cap K :t\in T\}$ is a partition of the set $K$.
By (\ref{EQ010}), for each $t\in T$ there is a bijection $\pi _t : [j _t]\cap J \rightarrow [j _t]\cap K$
and $\pi :=\bigcup _{t\in T}\pi _t :J \rightarrow K$ is a bijection.
If $j\in [j_t]$, then $\pi (j)\in [j_t]$ and, hence, $j\sim \pi (j)$. Thus (\ref{EQ027}) is true.

By the assumption, $\r$ is an equivalence relation on the set $I\cup J\cup K$ and $f\subset \r$.
So, since $\sim $ is the minimal equivalence relation on $I\cup J\cup K$ containing $f$
we have $\sim \;\subset \r$.
Thus, if $j\in J$, then by (\ref{EQ027}) we have $j\sim \pi (j)$ and, hence, $j\,\r \,\pi (j)$.
So, $\pi =\{ \la j ,\pi (j)\ra :j\in J\} \subset \r$ and the proof is finished.
\hfill $\Box$
\begin{ex}\label{EX003}\rm
{\it Lemma \ref{T042} is false, if $|I|\geq \o$}.
Let $I=\o$, $J=\{ j\}$ and $K=\{ k\}$, where
$j,k\not\in \o$, and let $f:I\cup J \rightarrow I\cup K$ be the bijection given by
$f=\{ \la j,0, \ra, \la 1,k\ra\} \cup \{ \la 2n,2n+2\ra :n\in \o\}\cup \{ \la 2n+1,2n-1\ra :n\geq 1 \}$.
Then $f_{rst}$ is the equivalence relation on the set $I\cup J\cup K$ determined by the partition
$\{ \{ j\} \cup \{2n:n\in \o\} , \{ k\} \cup \{2n+1:n\in \o\}\}$. Assuming that the lemma is true
we would have $\pi =\{ \la j,k\ra\}\subset f_{rst}$, that is $j \mathrel{f_{rst}} k$, which is false.
\end{ex}
\begin{prop}\label{T225}
If $(\ast )$ holds, $|I|< \o$ and $\Y$ and $\Z$ are $L_b$-structures disjoint from $\bigcup _{i\in I}X_i$, then

(a) If $\;\bcd _{i\in I}\X _i \du \Y \cong \bcd _{i\in I}\X _i \du \Z $, then $\Y \cong \Z $;

(b) $I(\bcd _{i\in I}\X _i \du \Y)\geq I(\Y )$.
\end{prop}
\dok
(a) Let $F:\bcd _{i\in I}\X _i \du \Y  \rightarrow \bcd _{i\in I}\X _i \du \Z $ be an isomorphism.
Let $\Y=\bcd _{j\in J}\Y _j$ and $\Z=\bcd _{k\in K}\Z _k$ be decompositions of the structures $\Y $ and $\Z $ into connectivity components,
where the index sets $I$, $J$ and $K$ are pairwise disjoint.
So, since $F:\bcd _{i\in I}\X _i \du \bcd _{j\in J}\Y _j \rightarrow \bcd _{i\in I}\X _i \du \bcd _{k\in K}\Z _k$ and $|I|<\o$,
by Fact \ref{T4015} we have $|J|=|K|$.
For $s\in I\cup J \cup K$ let us define
\begin{equation}\label{EQ241} \textstyle
\T_s=\left\{
           \begin{array}{rl}
           \X_s, & \mbox{ if } s\in I,\\
           \Y_s, & \mbox{ if } s\in J,\\
           \Z_s, & \mbox{ if } s\in K.\
           \end{array}
    \right.
\end{equation}
Then $F:\bcd _{s\in I\cup J}\T _s \rightarrow \bcd _{s\in I\cup K}\T _s$ is an isomorphism
and by Fact \ref{T4015} there is a bijection $f:I\cup J \rightarrow I\cup K$ such that
\begin{equation}\label{EQ242} \textstyle
\forall s\in I\cup J \;\;\T _s \cong \T _{f(s)}.
\end{equation}
Let $\r$ be the equivalence relation on the set $I\cup J \cup K$
given by $s\,\r \, t$ iff $\T _s \cong \T _t$.
By (\ref{EQ242}), for each $s\in I\cup J$ we have $s\,\r \, f(s)$,
and, hence, $f\subset \r$.
By Lemma \ref{T042} there is a bijection $\pi : J\rightarrow K$
such that for each $j\in J$ we have $j \,\r \, \pi (j)$,
that is $\T _j \cong \T _{\pi(j)}$
and, by (\ref{EQ241}), $\Y _j \cong \Z _{\pi(j)}$.
Thus, clearly, $\Y \cong \Z $.

(b) Let $|Y|\geq \o$, $\k :=I(\Y )$, and let $\Z _\a $, for $\a < \k $, be countable non-isomorphic models of $\Th (\Y )$.
By Fact \ref{T226}, for each $i\in I$ there is a countable connected structure $\X_i'\equiv \X _i$.
By  Proposition \ref{T200}(b) for each $\a <\k$ we have $\X _\a :=\bcd _{i\in I}\X _i ' \du \Z _\a \equiv \bcd _{i\in I}\X _i \du \Y$.
So, by (a) $\X _\a$, $\a <\k$, are countable non-isomorphic models of $\Th (\bcd _{i\in I}\X _i \du \Y)$.
\kdok
We remark that Claim (a) of  Proposition \ref{T225} is false if $|I|=\o$. For example, let $\X =\bcd _\o (1\du 2)$, $\Y =\bcd _\o 1$ and $\Z =\bcd _\o 2$.
An $L_b$-structure $\X$ is said to be {\it of diameter $\leq n$} iff for each $x,y\in X$ there is a path of length $\leq n$ from $x$ to $y$,
iff $\X\models \f _{\delta \leq n}$, where $\f _{\delta \leq n}:=\forall u,v\;\p _n(u,v)$.
We will say that {\it $\X$ is of finite diameter}, in notation $\delta (\X )<\o$, iff $\X$ is of diameter $\leq n$, for some $n\in \N$.
\begin{fac}\label{T210}
An $L_b$-structure $\X$ is of finite diameter iff each $\Y \models \Th (\X )$ is connected.
\end{fac}
\dok
If $\X \models \f _{\delta \leq n}$, then for each $\Y \models \Th (\X )$ we have $\Y \models \f _{\delta \leq n}$.

If $\X$ is not of finite diameter, then the type $\Phi (u,v):=\{ \neg \p _n(u,v):n\in \o\}$
is finitely satisfiable in $\X$ and, hence, there are $\Y \models \Th (\X )$ and $y,y'\in Y$
such that $\Y \models \Phi [y,y']$, which means that $y$ and $y'$ belong to different components of $\Y$,
and, hence, $\Y$ is not connected.
\hfill $\Box$
\begin{prop}\label{T201}
If $(\ast)$ holds and $\delta (\X_j )<\o <I(\X _j)$, for some $j\in I$, then $I(\bcd _{i\in I}\X _i)\geq I(\X _j)$.
\end{prop}
\dok
First suppose that $|X _i|\leq \o$, for $i\in I $.
Let $\X:=\bcd _{i\in I}\X _i$, $\k :=I(\X _j) >\o $
and let $\Y _\a$, for $\a <\k$, be non-isomorphic countable models of $\Th (\X _j)$, where $Y _\a \cap X=\emptyset$, for all $\a <\k$.
Since $\delta (\X_j )<\o$, by Fact \ref{T210} the structures $\Y _\a$, $\a <\k$, are connected.
By  Proposition \ref{T200}(b), $\X _\a := \bcd _{i\neq j}\X _i \du \Y _\a $, $\a <\k$, are countable models of $\Th ( \X)$.
Assuming that $I(\X )=\l <\k$ there would be $\o _1$-many isomorphic $\X _\a$'s (otherwise, we would have $\k \leq \l \o$)
say, using a re-enumeration, $\X _\a$, $\a <\o _1$.
So, for each $\a <\o _1$ we would have $\X _\a\cong \X _0$
and, by Fact \ref{T4015}, $\Y _\a$ would be isomorphic to some structure from the countable set $\{ \X _i :i\in I\setminus \{ j\}\}\cup \{ \Y _0\}$,
which is impossible.

If $|\X _i|>\o$, for some $i\in I$,
by Fact \ref{T226} there are pairwise disjoint connected structures $\X _i' \equiv \X _i$, such that $|X_i'|\leq \o$
and, clearly, $\delta (\X _j')<\o $.
By  Proposition \ref{T200}(b) we have $\bcd _{i\in I }\X _i' \equiv \X$,
since $I (\X _j')=I (\X _j)$ by (a) we have  $I(\bcd _{i\in I}\X _i')\geq I(\X _j')$,
that is $I(\bcd _{i\in I}\X _i)\geq I(\X _j)$.
\hfill $\Box$
\section{Closure under finite products and unions}\label{S3}
We regard the {\it binary direct product} $\times$ as a binary operation on the class $\Mod _{L_b}$ of $L_b$-structures
(as usual we identify $\X _0 \times \X _1$ and $\prod _{i<2}\X _i$)
and the {\it binary disjoint union} $\du$ as a partial binary operation on $\Mod _{L_b}$
($\X \du \Y$ is defined if $X\cap Y=\emptyset$).
If $\K \subset \Mod _{L_b}$ is a class of structures it is convenient to assume that it is $\cong$-invariant
(then for $\X,\Y \in \K$ there are $\X',\Y'\in \K$ such that $\X '\cong \X$, $\Y '\cong \Y$ and $X'\cap Y'=\emptyset$, so $\X ' \du \Y '$ is defined).
In Lemma \ref{T218} we describe the minimal closure of $\K$ under $\times$ and $\du$, denoted by $\la \K \ra$.
For convenience $\prod _{i<1}\X_0 :=\X_0$  is regarded as a finite product and $\bcd _{i<1}\X_0 :=\X_0$ as a finite union.
\begin{lem}\label{T218}
Let $n\in \N$ and $m_i\in \N$, for $i<n$. Then

(a) $\prod _{i<n} \bcd _{j<m_i}\X ^i_j = \bcd _{\bar j\in \prod _{i<n} m_i} \prod _{i<n}\X ^i_{j_i}$, whenever $\X ^i_j$, $j<m_i$, are pairwise disjoint $L_b$-structures, for each $i<n$;

(b) If the structures $\X ^i_j$ are connected, then $\X _{\bar j} := \prod _{i<n}\X ^i_{j_i}$, $\bar j \in \prod _{i<n} m_i$,
are the connectivity components of the structure from (a);

(c) If $\K \subset \Mod _{L_b}$ is closed under $\cong$, then $\la \K\ra $ is the class of all finite disjoint unions of finite products of structures from $\K$.
\end{lem}
\dok
(a) Let $\X ^i_j=\la X^i_j, \r^i_j\ra$, for $i<n$ and $j<m_i$, let $\X:=\prod _{i<n} \bcd _{j<m_i}\X ^i_j =\la X,\r\ra$ and $\Y:= \bcd _{\bar j\in \prod _{i<n} m_i} \prod _{i<n}\X ^i_{j_i} =\la Y,\s\ra$.
First we prove the equality of the domains: $\prod _{i<n} \bigcup _{j<m_i}X ^i_j=\bigcup _{\bar j \in \prod _{i<n} m_i} \prod _{i<n}X ^i_{j_i}$.
Since $X _{\bar j}:=\prod _{i<n}X ^i_{j_i} \subset X$, for all  $\bar j\in \prod _{i<n} m_i$, we have ``$\supset$".
If $\bar x=\la x_i :i<n\ra \in X$,
then for each $i<n$ there is a unique $j_i<m_i$ such that $x_i \in X^i_{j_i}$;
so $\bar j :=\la j_i :i<n \ra\in \prod _{i<n} m_i$ and $\bar x\in X_{\bar j}\subset Y$. Thus we have ``$\subset$".

Second we prove that $\r =\s$. For $\bar j \in \prod _{i<n} m_i$ let $\X _{\bar j} = \la \prod _{i<n}X ^i_{j_i}, \r_{\bar j}\ra $.
Since for a fixed $i$ the structures $\X ^i_j$, $j<m_i$, are pairwise disjoint,
the direct products $\X _{\bar j}$, $\bar j\in \prod _{i<n} m_i$, are pairwise disjoint,
which gives $\s =\bcd _{\bar j\in \prod _{i<n} m_i}\r_{\bar j}$.
So for $\bar x ,\bar y \in X$ we have:
$\la \bar x , \bar y\ra \in \r$
iff for each $i<n$ there is $j_i <m_i$ such that $\la x_i,y_i\ra \in \r ^i_{j_i}$;
iff there is $\bar j\in \prod _{i<n} m_i$ such that $\la x_i,y_i\ra \in \r ^i_{j_i}$, for each $i<n$;
iff there is $\bar j\in \prod _{i<n} m_i$ such that $\la \bar x , \bar y\ra \in \r_{\bar j}$;
iff $\la \bar x , \bar y\ra \in \s$.

(b) If the structures $\X ^i_j$ are connected the products $\X _{\bar j}$, $\bar j\in \prod _{i<n} m_i$, are connected (see e.g.\ Corollary 2.5 of \cite{KRP})
and pairwise disjoint.
Suppose that for different $\bar j,\bar k\in \prod _{i<n} m_i$ there are $\bar y \in X_{\bar j}$ and $\bar z \in X_{\bar k}$
such that $\bar y \mathrel{\r _{rst}} \bar z$
and let $i_0<n$, where $j_{i_0}\neq k_{i_0}$.
Thus there is a path from $\bar y $ to $\bar z$ in the product $\prod _{i<n} \bigcup _{j<m_i}\X ^i_j$
and it is easy to show (see e.g.\ Theorem 2.9 from \cite{KRP})
that there is $m\in \o$ such that for each $i<n$
there is a $m$-path from $y_i$ to $z_i$ in the factor $\bigcup _{j<m_i}\X ^i_j$.
In particular, there is a path from $y_{i_0}$ to $z_{i_0}$ in $\bigcup _{j<m_{i_0}}\X ^{i_0}_j$,
and, hence, there is $j_0<m_{i_0}$ such that $y_{i_0},z_{i_0}\in \X ^{i_0}_{j_0}$.
But, since  $\bar y \in X_{\bar j}$ and $\bar z \in X_{\bar k}$ we have $y_{i_0}\in \X ^{i_0}_{j_{i_0}}$ and $z_{i_0}\in \X ^{i_0}_{k_{i_0}}$,
which implies that $j_{i_0}=j_0=k_{i_0}$ and we have a contradiction.
Thus $\X _{\bar j} := \prod _{i<n}\X ^i_{j_i}$, $\bar j\in \prod _{i<n} m_i$, are maximal connected parts of $\X$.

(c) It is easy to show that $\la \K\ra =\bigcup _{n\in \o}\K_n$, where $\K _0 =\K$ and for $n\in \o$
$$
\K _{n+1}=\K _n \cup \{ \X _0 \times \X _1 \!: \X _0 , \X _1\in \K _n\}\cup \{\X _0 \du \X _1 : \X _0 , \X _1\in \K _n \land X_0 \cap X_1=\emptyset \}.
$$
If $\K ^*$ is the class of finite disjoint unions of finite products of structures from $\K$, the inclusion $\K ^* \subset \la \K\ra$ is evident.
For the converse by induction we prove that  $\K _n \subset \K ^*$, for each $n\in \o$.
For $n=0$ the claim is true, since for $\X_0 \in \K$ we have $\X _0=\prod _{i<1}\X_0 \in \K ^*$.

Suppose that  $\K _n \subset \K ^*$ and let $\X \in \K _{n+1}\setminus \K _n$ be constructed from $\X _0, \X _1\in \K _n$. Then by the induction hypothesis for each $i<2$ we have
$\X _i =\bcd _{j<m_i} \X ^i_j$, where $m_i \in \N$, $\X ^i _j=\prod _{k<n^i_j}\X ^i_{j,k}$ and $\X ^i_{j,k}\in \K$, for all $j<m_i$ and $k<n^i_j$.

If $\X=\X _0 \du \X _1$, then $X_0 \cap X_1=\emptyset$ and $\X =\bcd _{j<m_0} \X ^0_j \du \bcd _{j<m_1} \X ^1_j$ is a finite union of finite products of structures from $\K $; so $\X \in \K ^*$.

If $\X=\X _0 \times \X_1$,
then,  by (a),
$\X =\prod _{i<2} \bcd _{j<m_i}\X ^i_j=\bcd _{\bar j\in \prod _{i<2} m_i} \prod _{i<2}\X ^i_{j_i}$
$=\bcd _{\bar j\in \prod _{i<2} m_i} \prod _{i<2}\prod _{k<n^i_{j_i}}\X ^i_{{j_i},k}$
and since $\X ^i_{j,k}\in \K$, for all $i<2$, $j<m_i$ and $k<n^i_j$, we have $\X \in \K ^*$.
\kdok
We will say that a class $\K$ of structures is {\it good for products} iff $I(\prod _{i<n}\X _i)=\c$, whenever $n\in \N$, $\X_i \in \K$, for $i<n$, and $I(\X _i)=\c$, for some $i<n$.
\begin{te}\label{T219}
If $\K $ is an isomorphism-closed class of connected $L_b$-structures satisfying VC$^\sharp$
and $\K $ is good for products,
then each structure $\X \in \la \K \ra$ satisfies VC$^\sharp$.
More precisely, if $\X =\bcd _{i<n}\prod _{j<m_i}\X ^i_j$, then
\begin{equation}\label{EQ239}\textstyle
I(\X )=\left\{
                                                 \begin{array}{cl}
                                                         \leq 1, & \mbox{\rm if } I(\X ^i_j)\leq 1 ,\;\mbox{\rm for all } i<n \;\mbox{\rm  and }j<m_i ;\\[1mm]
                                                         \c ,    & \mbox{\rm if } I(\X ^i_j)=\c ,\;\mbox{\rm for some } i<n \;\mbox{\rm  and }j<m_i.
                                                 \end{array}
                                               \right.
\end{equation}
\end{te}
\dok
By Lemma \ref{T218}(c) each $\X \in \la \K \ra$ is of the form $\X =\bcd _{i<n} \X _i$, where $\X _i:=\prod _{j<m_i}\X ^i_j$, for $i<n$, and $\X ^i_j \in \K$, for $j<m_i$.
By Fact \ref{T226} there are connected models $\Y ^i_j \equiv \X ^i_j$ such that $|\Y ^i_j|\leq \o$, for all $i,j$;
thus, $I(\Y ^i_j)=I(\X ^i_j)\in \{ 0,1,\c\}$.
By Fact \ref{T043}(d) we have $\Y _i := \prod _{j<m_i}\Y ^i_j\equiv \X _i$, for $i<n$, and, hence, $I(\Y _i )=I(\X _i )$;
clearly, $\Y _i$'s are connected.
By  Proposition \ref{T200}(b) we have $\Y := \bcd _{i<n}\Y _i\equiv \X $ and, hence, $I(\Y )=I(\X )$.
If $I(\Y ^i_j)\leq 1$, for all $i<n$ and $j<m_i$,
then by Fact \ref{T043}(e) we have $I(\Y _i)\leq 1$, for all $i<n$,
and, by  Proposition \ref{T200}(c), $I(\Y )\leq 1$.
Otherwise, $I(\Y ^{i_0}_{j_0})=I(\X ^{i_0}_{j_0})=\c$, for some $i_0<n$ and $j_0<m_{i_0}$
and, by our assumption, $I(\X _{i_0})=\c$, that is $I(\Y _{i_0})=\c$.
Now, since $\Y _i$'s are connected, by  Proposition \ref{T225}(b) we have $I(\Y )=\c$
and, hence $I(\X )=\c$.
\hfill $\Box$
\section{Rooted  and initially finite trees}\label{S4}
We recall that a partial order $\X =\la X,\leq \ra$ is a {\it (model-theoretic) tree}
iff for each element $x$ of $X$ the set $(-\infty ,x]:=\{ y\in X : y\leq  x\}$
is linearly ordered by the relation $\leq $. The theory of trees is the $L_b$ theory
$\CT _{\rm tree}:=\{ \f _{\rm R},\f _{\rm AS},\f _{\rm T}\}\cup \{ \f _{\rm tree}\}$,
where $\f _{\rm tree}$ is the sentence $\f _{\rm tree}:=\forall v \;\forall u,w\;( u\leq v \land w\leq v \Rightarrow u\leq w \lor w\leq u)$.
If $\X =\la X,\leq  \ra $ is a tree, then $\parallel$ will denote the {\it comparability relation},
while $\not\perp$ will denote the {\it compatibility relation} on $X$ ($x\not\perp y$ iff there is $z$ such that $z\leq x$ and $z\leq y$),
definable in $\X$ by the formula $\f_{\not\perp }(u,v):=\exists w \; (w\leq u \land w\leq v )$.
\begin{fac}\label{T046}
If $\X =\la X,\leq \ra $ is a tree, then $\leq _{rs}\,=\;\,\parallel$ and $\leq _{rst}\;\,=\;\,\not\perp$.
\end{fac}
\dok
The equality $\leq _{rs}=\;\parallel$ is evident. Thus $\leq _{rst}$ is the transitive closure of $\parallel$.

If $x\not\perp y$, then $z\leq  x, y$, for some $z$, and, hence, $x\parallel z \parallel y$, so $x \leq _{rst} y$.

If $x \leq _{rst}\; y$, then there are $z_0,\dots,z_n $ such that $x=z_0\parallel z_1 \parallel\dots \parallel z_n =y$
and by induction on $n$ we show that $x\not\perp y$. For $n=1$ we have $x\parallel y$, so, $x\not\perp y$.
Suppose that the claim is true for $n$ and that $x=z_0\parallel z_1 \parallel\dots \parallel z_n \parallel z_{n+1}=y$.
Then by the hypothesis we have $x\not\perp z_n \parallel y$ and, hence, there is $z\leq x,z_n$.
If $z_n \leq y$, then $z\leq x,y$ and $x\not\perp y$. Otherwise we have $y\leq z_n$,
which implies that $z\parallel y$.
If $z\leq y$, then $z\leq x,y$ so, $x\not\perp y$.
If $y\leq z$, then $y\leq x$ and $x\not\perp y$ again.
\kdok
The equivalence classes $[x]$, $x\in X$, corresponding to the relation  $\leq _{rst}$
are the  (connected) components of $\X$.
If  $X/\!\! \leq _{rst}=\{ [x_i] :i\in I \}$, then
the suborders $\X _i :=\la [x_i] ,\leq \upharpoonright [x_i]\ra$, $i\in I$, are connected trees
and $\X$ can be regarded as the disjoint union of its components, $\X =\bcd _{i\in I}\X_i$.

If  $r\in X$ and $r\leq x$, for all $x\in X$,
then $r$ is the {\it root} of $\X$ and $\X$ is a {\it rooted tree}.
Then $\X ^+:=\X \setminus \{ r \}$ is a tree too and if $\X ^+$ has finitely many components, $\X$ will be called an {\it initially finite tree}.
\begin{fac}\label{T202}
The following properties are first-order properties in the class of  trees.

(a) $\X$ has $n$ components, for $n\in \N$;

(b) $\X$ is a disjoint union of linear orders;

(c) $\X$ is rooted and $\X ^+ $ has exactly $n$ components.
\end{fac}
\dok
(a) By Fact \ref{T046} the sentence $\f _n :=\exists v_0, \dots ,v_{n-1} \bigwedge _{0\leq i <j <n} v_i \perp v _j$, for $n\geq 2$,  says that $\X$ has at least $n$ components.
Thus the sentence $\neg \f _2$ says that $\X$ is connected and $\p _n:=\f _n \land \neg \f _{n+1}$
says that $\X$ has exactly $n$ components.

(b) By Fact \ref{T046}, the sentence $\f := \forall u,v \;( u\not\perp v \Rightarrow u \parallel v) $ says that components consist of comparable elements;
so, they are linear orders.

(c) Roughly, let $\f _{\rm loc.fin.}^n$ be the conjunction of a sentence saying that a root exists and the relativization of $\p _n $ to the set $X^+$.
\hfill $\Box$
\begin{ex}\rm \label{EX002}
{\it Connectedness is not a first-order property in the class of posets}.
Let $\X =\la Z ,\leq \ra$, where $Z$ is the set of integers
and $\leq$ the partial order on $Z$ given by $ \dots >-3<-2>-1<0>1<2>3 < \dots$.
Then $\X $ is a connected poset but it is not of finite diameter. By Fact \ref{T210} there is a disconnected model of $\Th (\X)$.
\end{ex}
\begin{fac}\label{T217}
A rooted tree $\X$ is $\o$-categorical iff the tree $\X ^+$ is $\o$-categorical.
\end{fac}
\dok
Let $r$ be the root of $\X$. If the tree $\X $ is $\o$-categorical, then, since $X=\{ r\} \cup X^+$ and for each $F\in \Aut (\X )$ we have $F(r)=r$ and $F[X^+]= X^+$,
by Proposition \ref{T216}(b) the tree $\X ^+ $ is $\o$-categorical.
Conversely, if the tree $\X ^+$ is $\o$-categorical, then,
since $\X$ is the lexicographical sum $\sum _2 \X_i$,
where $\X _0 =\la \{ r\}, \{\la r,r\ra\}\ra$ and $\X _1 =\X ^+$, the tree $\X$ is $\o$-categorical by  Proposition \ref{T200}(c).
\hfill $\Box$
\begin{prop}\label{T040}
If $\X _i $, $i\in I$, and $\Y _i $, $i\in I$, are rooted trees and $\prod _{i\in I}\X _i \cong \prod _{i\in I}\Y _i$,
then $\dot{\bigcup} _{i\in I}\X _i^+ \cong \dot{\bigcup} _{i\in I}\Y _i ^+$.
\end{prop}
\dok
For $i\in I$, let $\X _i =\la X_i ,\leq_{\X _i}\ra$, $\Y _i =\la Y_i ,\leq _{\Y _i}\ra$,
$0_{\X _i}:=\min \X _i$ and $0_{\Y _i}:=\min \Y_i$.
Let $X:=\prod _{i\in I}X _i$, $Y:=\prod _{i\in I}Y _i$,
$\X:=\prod _{i\in I}\X _i=\la X,\leq _\X\ra$ and $\Y:=\prod _{i\in I}\Y _i=\la Y,\leq _\Y\ra$.
For $i\in I$ let
$A_i:= \textstyle  \{ x \in X  : \forall j\in I \setminus \{ i \} \;x_j=0_{\X _j} \}$,
$B_i:=\textstyle  \{ y \in Y  : \forall j\in I \setminus \{ i \} \;y_j=0_{\Y _j} \}$
and let $A:=\bigcup _{i\in I}A_i$ and $B:=\bigcup _{i\in I}B_i$.
Let $\f (v)$ be a formula saying that the set $(-\infty ,v]$ is linearly ordered,
say $\forall u,w\;( u\leq v \land w\leq v \Rightarrow u\parallel w )$.
We show that
\begin{equation}\label{EQ238}
A=\{ x\in X : \X\models \f [x]\}\;\mbox{ and }\;  B=\{ y\in Y : \Y\models \f [y]\}.
\end{equation}
We prove only the first  equality.
Let $i\in I$ and $x\in A_i$.
If $x', x'' \leq_{\X } x$,
then for each $j\in I \setminus \{ i \}$ we have $x'_j, x''_j \leq_{\X _j } x_j=0_{\X _j}$
and, hence, $x'_j= x''_j =0_{\X _j}$.
Since $x_i', x_i''\leq_{\X _i} x_i$ and $\X _i$ is a tree, we have:
either $x_i'\leq _{\X _i} x_i''$, which gives $x'\leq _{\X }x''$;
or $x_i''\leq _{\X _i} x_i'$, which gives $x''\leq _{\X }x'$.
Thus $\X \models \f [x]$, for all $x\in A$.

If $x\not\in A$,
then there are different $i',i''\in I$
such that $0_{\X _{i'}}<_{\X _{i'}}x_{i'}$ and $0_{\X _{i''}}<_{\X _{i''}}x_{i''}$.
Let $x'\in X$, where $x_i'=0_{\X _i}$, for $i\neq i'$, and $x_{i'}'=x_{i'}$
and let $x''\in X$, where $x_i''=0_{\X _i}$, for $i\neq i''$, and $x_{i''}''=x_{i''}$.
Then we have $x',x''\leq_{\X } x$,
but $x'\not\leq_{\X } x''$ and $x''\not\leq_{\X } x'$;
so $\X \not\models \f [x]$ and the first equality is proved.

Let  $f:\X \rightarrow \Y$ be an isomorphism.
For $x\in X$ by (\ref{EQ238}) we have:
$x\in A$
iff $\X \models \f [x]$
iff $\Y \models \f [f(x)]$
iff $f(x)\in B$
iff $x\in f^{-1}[B]$.
So $A=f^{-1}[B]$
and, since $f$ is onto, $f[A]=B$.
Clearly, $0_\X :=\la 0_{\X_i}: i\in I\ra =\min \X$ and $0_\Y :=\la 0_{\Y_i}: i\in I\ra =\min \Y$
and, hence, $f(0_\X)=0_\Y$.
Also, for $i\in I$ we have $\A _i:=\la A_i ,\leq _\X \upharpoonright A_i\ra\cong \X _i$ and $\B _i:=\la B_i ,\leq _\Y \upharpoonright B_i\ra\cong \Y _i$,
which gives
\begin{equation}\label{EQ001'}
\A _i \setminus \{ 0_\X \}\cong \X _i \setminus \{ 0_{\X _i} \}=:\X _i ^+ \;\mbox{ and }\;
\B _i \setminus \{ 0_\Y \}\cong \Y _i \setminus \{ 0_{\Y _i} \}=:\Y _i ^+.
\end{equation}
Let $\A :=\la A ,\leq _\X \upharpoonright A\ra$ and $\B :=\la B ,\leq _\Y \upharpoonright B\ra$.
Since $f[A]=B$ and $f$ is an isomorphism,
the restriction $f\upharpoonright A :\A \to \B$ is an isomorphism
and, since $f(0_\X)=0_\Y$,
the restriction $f\upharpoonright (A\setminus \{ 0_\X \}) :\A\setminus \{ 0_\X \} \to \B\setminus \{ 0_\Y \}$ is an isomorphism too.
Thus
\begin{equation}\label{EQ002'}\textstyle
f\upharpoonright (A\setminus \{ 0_\X \}) :
\bigcup _{i\in I}\A _i\setminus \{ 0_\X \} \stackrel{\rm iso}{\longrightarrow} \bigcup _{i\in I}\B _i\setminus \{ 0_\Y \}.
\end{equation}
If $i\neq j\in I$, $x\in A _i\setminus \{ 0_\X \}$ and $x'\in A _j\setminus \{ 0_\X \}$,
then $x$ and $x'$ are $\leq _\X$-incomparable elements of $X$
so, by (\ref{EQ001'}) the suborder $\bigcup _{i\in I}\A _i\setminus \{ 0_\X \}$ of $\X$
is isomorphic to the disjoint union $\bcd _{i\in I} \X _i^+$
and, similarly,  $\bigcup _{i\in I}\B _i\setminus \{ 0_\Y \}\cong \bcd _{i\in I} \Y _i^+$,
which by (\ref{EQ002'}) gives $\bcd _{i\in I} \X _i^+ \cong \bcd _{i\in I} \Y _i^+$.
\hfill $\Box$
\begin{prop}\label{T041}
If $\X _r$, $r<n$, are initially finite trees, $\Y $ and $\Z $ are rooted trees and $(\prod _{r<n}\X _r )\times \Y \cong (\prod _{r<n}\X _r ) \times \Z $,
then $\Y \cong \Z $.
\end{prop}
\dok
W.l.o.g.\ we can assume that the structures $\X _r$, $r<n$, $\Y$ and $\Z$ have disjoint domains.
If $(\prod _{r<n}\X _r )\times \Y \cong (\prod _{r<n}\X _r ) \times \Z$,
then by Proposition \ref{T040} we have $\bcd _{r<n}\X _r ^+ \du \Y ^+ \cong \bcd _{r<n}\X _r ^+ \du  \Z ^+$.
Since the trees $\X _r$, $r<n$, are initially finite, the tree $\bcd _{r<n}\X _r ^+$ is a finite disjoint union of connected trees.
By  Proposition \ref{T225}(a) we have  $\Y ^+\cong \Z ^+$ and, clearly, $\Y \cong \Z $.
\hfill $\Box$
\begin{prop}\label{T211}
If $\X _r$, $r<n$, are countable rooted trees, then each of the following conditions implies that $I(\prod_{r<n}\X _r)=\c $.

(a) There is $s<n$ such that $I(\X_s)=\c$ and $\X_r$, for $r\neq s$, are initially finite;

(b) There is $s<n$ such that $I(\X_s)=\c$ and $\X_s$ is initially finite;

(c) There are $s<n$, $\Y\in \Mod (\Th (\X _s),\o)$ and a component $\C$ of $\Y ^+$ such that $I(\C)=\c$.
\end{prop}
\dok
(a) By Fact \ref{T043}(a), $\prod _{r< n}\X _r  \cong (\prod _{r\neq s}\X _r)\times \X _s$
and by Proposition \ref{T041}, if $\Y ,\Z\in \Mod (\Th (\X _s ),\o )$
and $\Y \not \cong \Z$,  then $(\prod _{r\neq s} \X _r ) \times \Y\not\cong (\prod _{r\neq s} \X _r ) \times \Z$.
So, by Fact \ref{T045} we have  $I(\prod _{r< n}\X _r )\geq I(\X _s)$.

(b) Let $\X _s ^+$ have $m$ components
and let $\Y _\a $, $\a <\c$, be non-isomorphic countable models of $\Th (\X _s )$.
W.l.o.g.\ we can assume that the structures $\X _r$, $r<n$, and $\Y_\a$, $\a <\c$,  have disjoint domains.
By  Fact \ref{T043}(d) $\Z _\a :=\prod_{r\in n\setminus \{s \}}\X _r\times \Y _\a$, $\a <\c$, are models of $\Th (\prod_{r<n}\X _r)$.
Suppose that $|\{\Z _\a  :\a <\c\}/\!\cong|=\k <\c$.
Then, since $\k\o<\c$, uncountably many of $\Z _\a$'s would be isomorphic;
using a re-enumeration, say $\Z _\a \cong \Z _0$, for $\a <\o _1$.
By Fact \ref{T202}(c), for $\a <\o_1$ we have $\Y _\a ^+=\dot{\bigcup} _{j<m}\Y_j ^\a$, where $\Y_j ^\a$, $j<m$, are the connectivity components of $\Y _\a ^+$.
For $r\in n\setminus \{s\}$ let $\X _r ^+=\bcd _{i\in I_r}\X ^r _i$ be the decomposition of the tree $\X _r ^+$ into connectivity components;
clearly, $|I_r|\leq \o$, for $r\in n\setminus \{s \}$.
Now for each $\a <\o _1$ we have $\Z _\a \cong \Z _0$;
so, by Proposition \ref{T040}
$\textstyle
\bcd _{r\in n\setminus \{s \}}\bcd _{i\in I_r}\X ^r _i \du \bcd _{j<m}\Y_j ^\a
\cong \bcd _{r\in n\setminus \{s\}}\bcd _{i\in I_r}\X ^r _i \du \bcd _{j<m}\Y_j ^0
$,
by Fact \ref{T4015} for each $\a <\o _1$ and $j<m$ we have $\tp (\Y ^\a _j)\in \{ \tp(\X ^r _i): r\in n\setminus \{s\} \land i\in I_r\}\cup \{ \tp (\Y_j ^0):j<m \}=:\{ \t _k :k <\k\}$,
and, clearly, $\k \leq \o$.
For each $\a <\o _1$ there is a unique non-decreasing function $c_\a :m \rightarrow \k$
such that $\tp (\Y _\a ^+)=\bcd _{j<m}\t _{c_\a (j)}$.
Since the structures $\Y _\a$ are non-isomorphic
their codes $c_\a$ are different.
But this is impossible because $\k ^m \leq \o$.

(c) Let $\Y ^+=\bcd _{j\in J}\Y_j \du \C$ be the decomposition of the tree $\Y  ^+$ into connectivity components,
let $\C _\a$, $\a <\c$, be non-isomorphic countable models of $\Th (\C)$
and let $\Y _\a ^+ :=\bcd _{j\in J}\Y_j \du \C _\a$, for $\a <\c$.
Then by  Proposition \ref{T200}(b) $\Y _\a ^+ \equiv\Y ^+$ and $\Y _\a :=(\Y _\a ^+)_r \equiv \Y $.
W.l.o.g.\ we can assume that the structures $\X _r$, $r<n$, and $\Y_\a$, $\a <\c$,  have disjoint domains.
By  Fact \ref{T043}(d) $\Z _\a :=\prod_{r\in n\setminus \{s \}}\X _r\times \Y _\a$, $\a <\c$, are models of $\Th (\prod_{r<n}\X _r)$.
Suppose that $|\{\Z _\a  :\a <\c\}/\!\cong|=\k <\c$.
Then, since $\k\o<\c$, uncountably many of $\Z _\a$'s would be isomorphic;
using a re-enumeration, say $\Z _\a \cong \Z _0$, for $\a <\o _1$.
Again, for $r\in n\setminus \{s\}$ let $\X _r ^+=\bcd _{i\in I_r}\X ^r _i$.
Now for each $\a <\o _1$ we have $\Z _\a \cong \Z _0$;
so, by Proposition \ref{T040}
$\textstyle
\bcd _{r\in n\setminus \{s \}}\bcd _{i\in I_r}\X ^r _i \du \bcd _{j\in J}\Y_j \du \C _\a
\cong \bcd _{r\in n\setminus \{s\}}\bcd _{i\in I_r}\X ^r _i \du \bcd _{j\in J}\Y_j \du \C _0
$.
By Fact\ref{T202}(a) the trees $\C _\a$ are connected,
by Fact \ref{T4015} for each $\a <\o _1$ we would have $\tp (\C _\a )\in \{ \tp (\X ^r _i ) :r\in n\setminus \{s \} \land i\in I_r\} \cup \{ \tp (\Y_j) : j\in J\} \cup \{ \tp (\C _0)\}$,
which is false because the trees $\C _\a$, $\a <\o _1$, are non-isomorphic, $|I_r|\leq \o$, for $r\in n\setminus \{s \}$, and $|J|\leq\o$.
Thus $I (\Th (\prod_{r<n}\X _r),\o)=\c $.
\kdok
Since initially finite trees are connected, by Theorem \ref{T219}, Fact \ref{T202}, Fact \ref{T043}(e) and Proposition \ref{T211} we have
\begin{te}\label{T224}
If $\K$ is the class of initially finite trees satisfying VC$^\sharp$,
then (\ref{EQ239}) holds and, hence, VC$^\sharp$ is true, for each poset $\X \in\la \K\ra$.
\end{te}
\section{VC$^\sharp$ in the closure of the class of rooted FMD trees}\label{S5}
We recall that, for a relational structure $\Y$ and $n\in \N$, $\f _\Y (n)$ denotes the number of non-isomorphic $n$-element substructures of $\Y$
and that the corresponding function $\f _\Y :\N \rightarrow \Card$ is called the {\it profile} of $\Y$.
$\Pa (\Y )$ denotes the set of all partial automorphisms (isomorphisms between substructures) of $\Y$ .
The structure $\Y$ is called
{\it monomorphic} iff $\f _\Y (n)=1$, for all $n\in \N$ such that $n\leq |Y|$;
{\it chainable} iff there is a linear order $<$ on $Y$ such that $\Pa (\la Y ,<\ra)\subset \Pa (\Y )$.
By the classical results of Fra\"{\i}ss\'{e} (see \cite{Fra}) and Pouzet \cite{Pou1}
an infinite relational structure $\Y$ is monomorphic iff $\Y$ is chainable iff $\Y$ is $\Pi _0$-definable
(i.e.\ definable by quantifier-free formulas) in some linear order $\la Y ,<\ra$.

We recall some results of Pouzet and Thi\'{e}ry from \cite{PT}.
If $\Y$ is a relational structure, then a partition $\{ Y_i:i\in I\}$ of its domain $Y$
is called a {\it monomorphic decomposition} of $\Y$ iff $\BK \cong \BH$, whenever $\BK$ and $\BH$ are finite substructures of $\Y$ such that
$|K \cap Y_i|=|H \cap Y_i|$, for all $i\in I$. There exists a minimal monomorphic decomposition of $\Y$
(see Propositions 1.6 and 2.12 of \cite{PT}) and if it is finite, then $\Y$ {\it admits a finite monomorphic decomposition}
(we write shortly: $\Y$ is an {\it FMD structure}).
By Theorem 1.8 of \cite{PT} the FMD structures are exactly the structures which are
$\Pi _0$-definable in the linear orders with unary predicates of the form $\X =\la X, < ,U_0^\X,\dots ,U_{m-1}^\X \ra$, where
$\{ U_0^\X,\dots ,U_{m-1}^\X \}$ is a partition of the domain $X$ consisting of convex sets.

By \cite{KFMD} if a complete theory $\CT$ of any relational language has an FMD model, then all models of $\CT$ are FMD;
this establishes the notion of an {\it FMD theory}. Moreover by Theorem 6.1 of \cite{KFMD} we have
\begin{fac}\label{T422}
{\rm VC}$^\sharp$ is true for each  FMD theory $\CT$.
\end{fac}
We note that the class of FMD structures is closed under substructures and $\Pi _0$-definability.
The class of FMD $L_b$-structures is closed under finite disjoint unions.
\begin{prop}\label{T049}
Let $\Y$ be a poset. Then we have

(a) $\Y$ is chainable  $\Leftrightarrow$ $\Y$ is monomorphic $\Leftrightarrow$ $\Y$ is a chain or an antichain;

(b) $\Y$ is an FMD poset iff $\Y=\sum _{\BI}\Y _i$, where $\BI$ is a finite poset and $\Y _i$ is a chain or an antichain, for each $i\in I$;

(c) $\Y$ is an FMD tree iff $\Y=\sum _{\BI}\Y _i$, where $\BI$ is a finite tree, $\Y _i$ is a chain for each non-maximal element $i\in I$,
and $\Y _i$ is a chain or an antichain, for each maximal element $i\in I$;

(d) An FMD tree $\Y=\sum _{\BI}\Y _i$ is connected iff the tree $\BI$ has a root $i_0$ and $\Y _{i_0}$ is a linear order (i.e.\ $\Y _{i_0}$ is not an antichain of size $>1$);

(e) An FMD tree $\Y=\sum _{\BI}\Y _i$ is rooted iff the tree $\BI$ has a root $i_0$ and $\Y _{i_0}$ is a linear order having a minimum.
\end{prop}
\dok
(a) Clearly, chainable structures are monomorphic.
Let $\Y$ be monomorphic and $x,y\in Y$, where $x\neq y$.
If $x \parallel y$
then (since $\Y$ is 2-monomorphic) each two elements of $Y$ are comparable and, hence, $\Y$ is a chain.
Otherwise, each two elements of $Y$ are incomparable and $\Y$ is an antichain.
Finally, each chain is chainable and
each antichain $\Y$ is chainable by any linear order on $\Y$,
since $\Pa (\Y )$ is the set of all bijections between subsets of $Y$.

(b) By  Corollary 2.26 of \cite{PT},
an $L_b$-structure $\Y$ is an FMD structure
iff $\Y=\sum _{\BI }\Y_i$, where $\BI$ is a finite $L_b$-structure
and $\Y_i$, $i\in I$, are chainable $L_b$-structures.
So, $\Y$ is an FMD poset
iff $\Y=\sum _{\BI}\Y _i$, where $\BI=\la I, \r\ra$ is a finite $L_b$-structure, say reflexive,
and $\Y _i$, $i\in I$, are chainable $L_b$-structures. By (a) and since $\Y _i$'s are suborders of $\Y$, they are chains or antichains.
Taking $y_i\in Y_i$, for $i\in I$, it is easy to check that $i\mapsto y_i$ is an isomorphism from $\BI$ to $\Y \upharpoonright \{ y_i :i\in I\}$,
which implies that $\BI$ is a poset.

(c) If $\Y$ is an FMD tree,
then by (b) $\Y=\sum _{\BI}\Y _i$,
where $\BI$ is a finite poset and $\Y _i$ is a chain or an antichain, for each $i\in I$.
Since $\BI$ embeds in $\Y$, $\BI$ must be a tree.
Assuming that $i<j $ and that $Y_i$ is an antichain of size $>1$,
for an $x\in Y_j$ $(-\infty ,x]$ would not be a chain.
Conversely if $\Y$ is of the specified form, by (b) it is a FMD poset, and an easy discussion shows that $(-\infty ,x]$ is a chain, for each $x\in X$.

(d) Let us choose $y_i\in Y_i$, for $i\in I$.
Let $\Y=\sum _{\BI}\Y _i$ be a connected FMD tree.
Since $\BI$ is a finite tree it has a minimal element, say $i_0$.
For $i\in I\setminus \{ i_0\}$ by Fact \ref{T046} we have $y_i \not\perp _{\Y} y_{i_0}$
and, hence, there is $i_1\in I$ and $z\in Y _{i_1}$ such that $z\leq _\Y y_i , y_{i_0}$.
Assuming that $i_1 \neq i_0$ we would have $y_{i_1}\leq _\Y y_i , y_{i_0}$
and, since $Y_{i_1}\cap Y_{i_0}=\emptyset$, $y_{i_1}< _\Y y_{i_0} $, we would have  $i_1< _\BI i_0 $,
which is impossible by the minimality of $i_0$.
Thus $i_1 =i_0$, which gives $y_{i_0}\leq _\Y y_i $, in fact, $y_{i_0}< _\Y y_i $, and, hence, $i_0< _\BI i $.
So, $i_0$ is the root of $\BI$.
Assuming that $\Y _{i_0}$ is an antichain of size $>1$,
by (a) $i_0$ would be a maximal element of $\BI$
and, hence $I=\{ i_0\}$ and $\Y =\Y _{i_0}$, which is impossible.
Thus $\Y _{i_0}$ is either a singleton or a linear order of size $>1$.

Conversely, let $i_0$ be the root of $\BI$,
and let $\Y _{i_0}$ be a linear order.
By Fact \ref{T046}, for $x,y\in Y$ we have to prove that $x \not\perp _{\Y} y$.
If $x,y\in \Y _{i_0}$, then $x \parallel y$ and, hence, $x \not\perp _{\Y} y$
If $x\in \Y _{i}$ and $y\in \Y _{j}$, where $i,j \neq i_0$,
then, since $i_0 <_{\BI}i,j$ we have $y_{i_0} <_\Y x,y$ and  $x \not\perp _{\Y} y$ again
If $x\in \Y _{i_0}$ and $y\in \Y _{i}$, where $i \neq i_0$,
then, $x <_\Y y$ and  $x \not\perp _{\Y} y$.
(e) follows from (d).
\hfill $\Box$
\begin{prop}\label{T050}
If $\X _r$, $r<n$, are rooted FMD trees and $I(\X _s )= \c$, for some $s< n$, then $I(\prod _{r<n}\X _r)= \c$.
\end{prop}
\dok
It is sufficient to prove that
{\it for each rooted FMD tree $\Y$ satisfying $I(\Y )= \c$ there is a component $\Z$ of $\Y ^+$ such that $I(\Z )= \c$}.\footnote{This is false without FMD;
see Theorem \ref{T203}, add a root to a disjoint union of countably many countable non-isomorphic $\o$-categorical linear orders.}
Then there is a component $\Z$ of $\X _s ^+$ such that $I(\Z)= \c$ and, by Proposition \ref{T211}(c),  $I(\prod _{r<n}\X _r)= \c$.

So, by (c) and (e) of Proposition \ref{T049} we have $\Y=\sum _{\BI}\Y _i$, where $\BI$ is a finite tree,
with a root $i_0$ and $\Y _{i_0}$ is a linear order having a minimum, say $r:=\min Y_{i_0}$.
Then $r$ is the root of $\Y$ and $\Y ^+ =\Y \setminus \{ r\}$. Since $I(\Y )= \c$, by Facts \ref{T217} and  \ref{T422} we have $I(\Y ^+)= \c$.
If $Y_{i_0}' :=Y_{i_0} \setminus \{ r\} \neq \emptyset$,
then defining $\Y _i' :=\Y _i$, for $i\in I \setminus \{ i_0\}$,
we have $\Y ^+ =\sum _{\BI}\Y _i'$
and, since $Y_{i_0}'$ is a linear order,
by Proposition \ref{T049}(d) the tree $\Z :=\Y ^+$ is connected; so, we are done.

Otherwise we have $Y_{i_0} = \{ r\}$ and, hence, $I\setminus \{ i_0\}\neq \emptyset$ and $\Y ^+=\sum _{i>_{\BI}i_0}\Y _i$.
Since $I$ is a finite set $\BI$ is a set-theoretic tree;
so, $|\Lev _1 (\BI)|=m\in \N$.
Let $J:=\{ i\in \Lev _1 (\BI) : (\Y _i \mbox{ is not a chain}) \mbox{ or } (i \mbox{ is a max.\ element of } \BI \mbox{ and } |Y_i|=1 \}$.
Assuming that $J=\Lev _1 (\BI)$,
by Proposition \ref{T049}(c) each $i\in \Lev _1 (\BI)$ would be a maximal element of $\BI$
and, hence, $I=\{ i_0\}\cup \Lev _1 (\BI)$.
In addition we would have $\Y \cong \o ^{<2}$ (rooted infinite antichain);
but it is evident that the tree $\o ^{<2}$ is $\o$-categorical
(while $I(\Y )= \c$). Thus $J \varsubsetneq \Lev _1 (\BI)$. Let us take an enumeration
\begin{equation}\label{EQ220}
\Lev _1 (\BI)=\{ i_1, \dots ,i_l ,i_{l+1}, \dots , i_m\}, \mbox{ where } J=\{ i_{l+1}, \dots , i_m\} .
\end{equation}
So $l\geq 1$ and
\begin{equation}\label{EQ221}
\forall k\in \{ 1,\dots ,l\} \;\;\Big( (\Y _{i_k} \mbox{is a chain}) \land (|Y _{i_k}|>1 \lor  \exists j\in I \; j>i_k)\Big).
\end{equation}
Let $I_k:=\{ i\in I :i_k \leq _{\BI} i\}$, for  $k\in \{ 1,\dots ,l\}$.
Clearly $\BI _k :=\la I_k, \leq _{\BI} \upharpoonright I_k\ra$ is a tree with root $i_k$
and, by Proposition \ref{T049}(d) and (\ref{EQ221}), $\Z _k := \sum _{\BI _k}\Y _i$ is a connected FMD tree.
Since $\Lev _1 (\BI)$ is an antichain in $\BI \setminus \{ i_0\}$, $\Z _k$ is a component of $\Y ^+$.

By (\ref{EQ220}) for  $k\in \{ l+1,\dots ,m\}$ $i_k$ is a maximal element of $\BI$,
which gives $\Y ^+ =\bigcup _{1\leq k \leq l} \Z _k \cup \bigcup _{l+1 \leq k \leq m}\Y _{i_k}$
and since $\Y _{i_k}$, for $k\in \{ l+1,\dots ,m\}$, are antichains (or singletons)
the union $\bigcup _{l+1 \leq k \leq m}\Y _{i_k}$ is isomorphic to one antichain, $\A $.
So the tree $\Y ^+$ can be presented as the disjoint union $\Y ^+ =\bcd _{k\in \{ 1,\dots ,l\}} \Z _k \du \A$.
Suppose that all the trees $\Z _k$, $k\in \{ 1,\dots ,l\}$, are $\o$-categorical.
Then, since the antichain $\A $ is $\o$-categorical too,
by  Proposition \ref{T200}(c) the tree $\Y ^+$ would be $\o$-categorical,
which is false since $I(\Y ^+)= \c$.
Thus there is $k\in \{ 1,\dots ,l\}$ such that the tree $\Z _k$ is not $\o$-categorical,
and, since it is an FMD tree, by Fact \ref{T422} we have $I(\Z _k )= \c$.
\kdok
\noindent
Since rooted trees are connected, by Fact \ref{T422},  Proposition \ref{T050} and Theorem \ref{T219} we have
\begin{te}\label{T221}
If $\K$ is the class of rooted FMD trees,
then (\ref{EQ239}) holds and, hence, VC$^\sharp$ is true, for each poset $\X \in\la \K\ra$.
\end{te}
\section{VC$^\sharp$ in the closure of the class of linear orders}\label{S6}
By Theorem \ref{T224} VC$^\sharp$ is true for each poset from $\la \LO _0\ra$, where $\LO _0$ is the class of linear orders with a smallest element.
Here we prove the same for the class $\LO $.
\begin{prop}\label{T038}
If $\X _i $ and $\Y _i $, for $i\in I$, are linear orders, then
$$\textstyle
\prod _{i\in I}\X _i \cong \prod _{i\in I}\Y _i \;\;\Leftrightarrow \;\; \exists \pi \in \Sym (I) \;\;\forall i\in I \;\;\X _i \cong \Y_{\pi (i)}.
$$
\end{prop}
\dok
Let $\X _i =\la X_i ,\leq_{\X _i}\ra$ and $\Y _i =\la Y_i ,\leq _{\Y _i}\ra$, for $i\in I$;   $X:=\prod _{i\in I}X _i$, $Y:=\prod _{i\in I}Y _i$,
$\X:=\prod _{i\in I}\X _i=\la X,\leq _\X\ra$ and $\Y:=\prod _{i\in I}\Y _i=\la Y,\leq _\Y\ra$.
Let  $f:\X \rightarrow \Y$ be an isomorphism, $a=\la a_i :i\in I\ra \in X$ and $f(a)=b=\la b_i :i\in I\ra$. We show that
\begin{equation}\label{EQ240}\textstyle
f[\prod _{i\in I} [a_i,\infty)]=\prod _{i\in I} [b_i,\infty) \;\;\mbox{ and }\;\;f[\prod _{i\in I} (-\infty ,a_i]]=\prod _{i\in I} (-\infty ,b_i].
\end{equation}
If $x=\la x_i:i\in I\ra\in \prod _{i\in I} [a_i,\infty)$,
then $a_i \leq _{\X _i} x_i$, for all $i\in I$,
and, hence, $a\leq _{\X }x$,
which gives $b=f(a)\leq _{\Y }f(x)=:y=\la y_i :i\in I\ra$.
Thus $b_i \leq _{\Y _i}y_i$, for all $i\in I$;
thus, $f(x)=y \in \prod _{i\in I} [b_i,\infty)$.
Conversely, if $y=\la y_i :i\in I\ra\in \prod _{i\in I} [b_i,\infty)$,
then $b_i \leq _{\Y _i}y_i$, for all $i\in I$,
and, hence, $b\leq _{\Y } y$.
Since $f$ is an surjection there is $x=\la x_i:i\in I\ra\in X$ such that $f(x)=y$
and, since $f$ is an isomorphism and $f(a)=b\leq _{\Y }y=f(x)$,
we have $a\leq _{\X }x$.
Thus $a_i \leq _{\X _i} x_i$, for all $i\in I$,
and, hence, $x\in \prod _{i\in I} [a_i,\infty)$,
which gives $y= f(x)\in f[ \prod _{i\in I} [a_i,\infty)]$.
Thus the first equality in (\ref{EQ240}) is proved and the proof of the second is dual.

By (\ref{EQ240}) the restrictions $f\upharpoonright \prod _{i\in I} [a_i,\infty): \prod _{i\in I} [a_i,\infty)\rightarrow \prod _{i\in I} [b_i,\infty)$ and
$f\upharpoonright \prod _{i\in I} (-\infty ,a_i]: \prod _{i\in I} (-\infty ,a_i]\rightarrow \prod _{i\in I} (-\infty ,b_i]$
are isomorphism of products of linear orders with a smallest element.
For $i\in I$ let
\begin{eqnarray*}
A_i & := & \textstyle  \{ x \in X  : a_i\leq _{\X _i}x_i \,\land \,\forall j\in I \setminus \{ i \} \;x_j=a_j \},\\
B_i & := & \textstyle  \{\, y \in Y  : b_i\,\leq _{\Y _i}y_i \;\land \;\forall j\in I \setminus \{ i \} \;y_j=b_j \},\\
C_i & := & \textstyle  \{ x \in X  : x_i\leq _{\X _i}a_i \,\land \,\forall j\in I \setminus \{ i \} \;x_j=a_j \},\\
D_i & := & \textstyle  \{\, y \in Y  : y_i\,\leq _{\Y _i}b_i \;\land \;\forall j\in I \setminus \{ i \} \;y_j=b_j \}.
\end{eqnarray*}
Using the proof of Proposition \ref{T040} we show that there are permutations $\pi,\sigma \in \Sym (I)$ such that
\begin{equation}\label{EQ003}
\forall i\in I \; ( f[A_i]=B_{\pi (i)}  \;\mbox{ and }\; f[C_i]= D_{\sigma (i)}).
\end{equation}
By (\ref{EQ240}) $f\upharpoonright \prod _{i\in I} [a_i,\infty): \prod _{i\in I} [a_i,\infty)\rightarrow \prod _{i\in I} [b_i,\infty)$ is an isomorphism;
clearly, $\prod _{i\in I} [a_i,\infty)$ and $\prod _{i\in I} [b_i,\infty)$ are products of rooted trees and $a$ and $b$ are their smallest elements.
If $i\neq j\in I$, $x\in A _i\setminus \{ a \}$ and $x'\in A _j\setminus \{ a \}$,
then $x$ and $x'$ are incomparable elements of $\prod _{i\in I} [a_i,\infty)$ .
So, since $\A _i \cong [a_i,\infty)$ and $\B _i \cong [b_i,\infty)$,
the suborders $\bigcup _{i\in I}(\A _i\setminus \{ a \})\subset \X$ and $\bigcup _{i\in I}(\B _i\setminus \{ b \})\subset \Y$ are disjoint unions of linear orders.
By (\ref{EQ002'}) (see the proof of Prop.\ \ref{T040})
we have $f\upharpoonright (\bigcup _{i\in I}A _i\setminus \{ a \}) : \bcd _{i\in I}\A _i\setminus \{ a \} \rightarrow \bcd _{i\in I}\B _i\setminus \{ b \}$ is an isomorphism,
by Fact \ref{T4015} there is $\pi \in \Sym (I)$ such that for each $i\in I$ we have $f[ A _i\setminus \{ a \}]=B _{\pi (i)}\setminus \{ b \}$
and, since $f(a)=b$, $f[A _i]=B _{\pi (i)}$.
The proof of the existence of $\sigma \in \Sym (I)$ satisfying (\ref{EQ003}) is dual.

For $i\in I$ the suborder $C_i \cup A_i= \{ x \in X  : \forall j\in I \setminus \{ i \} \;x_j=a_j \}$ of $\X$
is a linear order isomorphic to $\X_i$
and by (\ref{EQ003}),  $f[C_i \cup A_i]=D_{\sigma (i)}\cup B_{\pi (i)}$,
which implies that $D_{\sigma (i)}\cup B_{\pi (i)}$ is a linear suborder of $\Y$ isomorphic to $\X_i$.
Assuming that $\pi (i)\neq \sigma (i)$
we would have that $D_{\sigma (i)}\cup B_{\pi (i)}$ is not a linear order.
Thus $\pi (i)= \sigma (i)$, for all $i\in I$,
and by (\ref{EQ003}), $\X _i \cong f[C_i \cup A_i]=D_{\pi (i)}\cup B_{\pi (i)}\cong \Y _{\pi (i)}$, for all $i\in I$.
\hfill $\Box$
\begin{prop}\label{T047}
If $\,n\in \N$, $\X _i$, $i< n$, are countable linear orders and $I(\X _j )=\c$, for some $j<n$,
then $I(\prod _{i< n}\X _i)=\c$.
\end{prop}
\dok
By Fact \ref{T043}(a), $\prod _{i< n}\X _i  \cong (\prod _{i\neq j}\X _i)\times \X _j$
and by Fact \ref{T045} we will have  $I(\prod _{i< n}\X _i )\geq I(\X _j )$, if for all $\Y ,\Y '\in \Mod (\Th (\X _j ),\o )$
from $(\prod _{i\neq j} \X _i ) \times \Y\cong (\prod _{i\neq j} \X _i ) \times \Y'$ it follows that $\Y \cong \Y '$.
So, let $(\prod _{i\neq j} \X _i ) \times \Y\cong (\prod _{i\neq j} \X _i ) \times \Y'$,
let $I:=n\setminus \{ j\}$, $J:=\{j\}$, $K:=\{ n\}$ and for $s\in I\cup J \cup K$ let
\begin{equation}\label{EQ020}\textstyle
\Z_s=\left\{
           \begin{array}{ll}
           \X_s, & \mbox{ if } s\in I,\\
           \Y, & \mbox{ if } s=j,\\
           \Y', & \mbox{ if } s=n.\
           \end{array}
    \right.
\end{equation}
Then $\prod _{s\in I\cup J}\Z _s= (\prod _{i\neq j} \X _i ) \times \Y\cong (\prod _{i\neq j} \X _i ) \times \Y'=\prod _{s\in I\cup K}\Z _s$
and, by Proposition \ref{T038}, there is a bijection $f:I\cup J \rightarrow I\cup K$ such that $\Z _s\cong \Z _{f(s)}$, for all $s\in I\cup J$.
Let $\r$ be the equivalence relation on the set $I\cup J \cup K$ defined by $s\,\r \, t$ iff $\Z _s \cong \Z _t$.
Then we have $s\,\r \, f(s)$, for all $s\in I\cup J$
and by Lemma \ref{T042} there is a bijection $\pi : J\rightarrow K$ such that for each $j\in J$ we have
$j \,\r \, \pi (j)$.
Clearly, $\pi (j)=n$ and, hence, $\Z _j \cong \Z _n$
and, by (\ref{EQ020}), $\Y  \cong \Y '$.
\kdok
Since linear orders are connected, by Proposition \ref{T047} and  Theorem \ref{T219} we have
\begin{te}\label{T220}
If $\K$ is the class of linear orders,
then (\ref{EQ239}) holds and, hence, VC$^\sharp$ is true, for each poset $\X \in\la \K\ra$.
\end{te}
\section{VC$^\sharp$ for disjoint unions of linear orders}\label{S7}
First we recall some basic facts concerning $\o$-categorical linear orders.
For $n\in \N$ let $\{ D_i :i<n\}$ be a partition of the rational line $\Q =\la Q, <\ra$ into dense subsets,
let $\t _i$, $i < n$,  be linear order types
and let $\X _q$, for $q\in Q$, be linear orders with pairwise disjoint domains
such that $\tp (\X _q )=\t _i$ iff $q\in D_i$.
The linear order $\sum _{q\in \Q} \X  _q$ is the {\it shuffle} of the set
$\{\t _i :i < n\}$ and it is denoted by $\s (\{\t _i :i < n\})$.
The notation $\s (\{\t _i :i < n\})$ is legal;
i.e., it is (up to isomorphism) independent of the choice of the partition $\{ D_i :i<n\}$,
of the enumeration of $\t _i$'s, and of the choice of their representatives, $\X _q$'s (see \cite{Rosen}, p.\ 116).

By the characterization of Rosenstein (see \cite{Rosen}, p.\ 299),
the class of at most countable $\o$-categorical linear order types
is the smallest class $\CC$ of order types containing:
1;
$\t _1 +\t _2$, whenever $\t _1,\t _2 \in \CC$; and
$\s (F)$, for each finite $F \subset \CC$.
The {\it rank} of $\t \in \CC$ is defined by $\rank (\t):=\min \{ n\in \o : \t \in \CC _n\}$,
where $\CC_0 =\{ 1\}$ and $\CC _{n+1}=\{ \t _1 +\t _2 : \t _1,\t _2 \in \CC_n \} \cup \{ \s (F) : F\in [\CC_n]^{<\o}\}$.
The collection of {\it predecessors} $\pred (\t)$ of a type $\t \in \CC$ is defined in the following way:
$\pred (1)=\emptyset$ and if $\rank (\t )=n+1$, then
\begin{eqnarray*}
\pred (\t ) & := & \textstyle \bigcup _{\t =\t _1 +\t _2 \land \rank (\t _1),\rank (\t _2 )\leq n} \{ \t _1 ,\t _2\}\cup  \pred (\t _1 )\cup  \pred (\t _2 )\cup \\
            &    & \textstyle \bigcup _{\t =\s (F) \land \forall \s \in F \;\rank (\s )\leq n} (F \cup \bigcup _{\s \in F}  \pred (\s )).
\end{eqnarray*}
As in \cite{Rosen} the  convex parts of a linear order $\X$ are called the {\it intervals} of $\X$.
\begin{fac}\label{T206}
(a) The collections $\CC _n$ and $\{ \t \in \CC : \rank (\t)=n\}$, $n\in \o$, are finite;

(b) If $\t =\tp (\X )\in \CC$, $\rank (\t )= n$ and $\Y$ is an interval of $\X$, then $\s :=\tp (\Y )\in \CC$ and $\rank (\s )\leq 2n+1$ (see \cite{Rosen}, p.\ 299);

(c) If $\s \in \pred (\t)$, then $\rank (\s )< \rank (\t )$, $\pred (\s)\subset \pred (\t)$ and  $\s$ is isomorphic to a proper interval of $\t$;

(d) For each infinite collection $\CF \subset \CC$ there are order types $\s _n \in \CC$, for $n\in \o$, such that for each $n\in \o$ we have
$\rank (\s _n)=n$, $\s _n\in \pred (\s _{n+1})$ and
$\s _n\in \pred (\t )$, for some $\t \in \CF$
 (see \cite{Rosen}, p.\ 322).
\end{fac}
\begin{fac}\label{T205}
(a) The theory of each $\o$-categorical linear order is finitely axiomatizable (L \"{a}uchli and Rosenstein \cite{Rosen1}).

(b) The theory of an $\o$-categorical tree is finitely axiomatizable iff the tree is finite-branching (Schmerl  \cite{Sch0}, Theorem 2.2).
\end{fac}
\begin{te}\label{T203}
If an infinite structure $\X $ is a disjoint union of linear orders, then $I(\X )\in \{ 1,\c \}$.
The theory $\Th (\X )$ is $\o$-categorical iff all the components of $\X$ are $\o$-categorical or finite and the collection of their order types is finite.
Then the theory $\Th (\X )$ is finitely axiomatizable iff $\X$ has finitely many components.
\end{te}
\dok
By the L\"{o}wenheim-Skolem theorem we can assume that the structure $\X $ is countable.
By Fact \ref{T202}(c), $\X =\bcd _{i\in I}\X_i$, where $\X _i$, for $i\in I$, are linear orders.
If $I(\X _j)=\c$, for some $j\in I$, then by  Proposition \ref{T201} and since $\X _j$ is of finite diameter we have $I(\X )=\c$.

Otherwise, by Rubin's theorem, the linear orders $\X _i$, $i\in I$, are finite or $\o$-categorical and for $\t _i :=\tp (\X _i)$ we have $\t _i\in \CC$, for all $i\in I$.

{\it Case 1}: $|\{ \t _i :i\in I\}|<\o$. Then the theory  $\Th (\X )$ is $\o$-categorical by  Proposition \ref{T240}.
Clearly, the set of the components of the tree $\X$ is finite iff it is finite-branching iff, by Fact \ref{T205}(b), $\Th (\X)$ is finitely axiomatizable.

{\it Case 2}:  $|\{ \t _i :i\in I\}|=\o$.
Then by Fact \ref{T206}(d) there are order types $\s _n \in \CC$, for $n\in \o$, such that for each $n\in \o$ we have:
(i) $\rank (\s _n)=n$,
(ii) $\s _n\in \pred (\s _{n+1})$ and
(iii) $\s _n\in \pred (\t _{i_n})$, for some $i_n \in I$.

For each $n\in \o$ we choose a component $\X _{i_n}$ of $\X$ such that $\tp (\X _{i_n})=\t _{i_n}$.
By (iii) for each $n\in \o$ the component $\X _{i_n}$ contains an interval $J_n$
such that $\tp (J_n)=\s _n$.
By (ii) and Fact \ref{T206}(c) $J_n$ is isomorphic to a proper interval of $J_{n+1}$;
thus, for each $n\in \o$ there is an embedding $f_n :J_n \rightarrow J_{n+1}$
such that $f_n [J_n]$ is a proper interval of $J_{n+1}$.
It is evident that the limit of that direct system (see \cite{CK}, p.\ 243)
is isomorphic to the linear order $J' =\bigcup _{n\in \o}J_n'$,
where $J_0' \subset J_1 ' \subset \dots $ is a strictly increasing family of intervals of $J'$
and $\tp (J_n')=\s _n$.

{\it Subcase 2.1}: $\forall n\in \o \; \exists a,b\in J' \; a< J_n' <b$.
Assuming that there is $c=\min J'$
we would have $c\in J_n '$, for some $n\in \o$, which is impossible.
Similarly, $\max J'$ does not exist.
Let $\la a_k :k\in \o\ra$  and $\la b_k :k\in \o\ra$
be a coinitial and a cofinal sequence in $J'$ respectively,
where $\dots <a_1 <a_0 <b_0 <b_1<\dots $.
Then $J' =\bigcup _{k\in \o}[a_k ,b_k]$.
For $k\in \o$ there are $n',n''\in \o$ such that $a_k \in J_{n'}'$ and $b_k \in J_{n''}'$;
thus for $n:=\max \{ n',n''\}$ we have $a_k,b_k\in J_n'$,
and by the convexity of $J_n'$ we have $[a_k,b_k]\subset J_n'$; so
\begin{equation}\label{EQ205}
\forall k\in \o \;\;\exists n\in \o \;\; ([a_k,b_k] \mbox{ is an interval in } J_n') .
\end{equation}
By Fact \ref{T206}(b) we have $\tp([a_k,b_k])\in \CC$
and, by Fact \ref{T205}(a), there is a sentence $\f _k$ axiomatizing the theory $\Th ([a_k,b_k])$,
w.l.o.g.\ we assume $\Var (\f _k)\subset \{ v_i:i\in \o\}$.
Let $L:=\{ \leq \} \cup \{ c_k :k\in \o\} \cup \{ d_k :k\in \o\}$, where $c_k$ and $d_k$ are constants.
For $k\in \o$ let $\theta _k (u):=c_k \leq u \leq d_k$,
let $\f _k^{\theta _k}$ be the corresponding relativization
and $\CT := \Th (\X )\cup \{ c_{k+1}<c_k<d_k<d_{k+1} :k\in \o\} \cup \{ \f _k^{\theta _k} :k\in \o\}$.
We show that the theory $\CT$ is finitely satisfiable.
Let $\Sigma \in [\Th (\X )]^{<\o}$, $m\in \o$ and
$\Sigma ' := \Sigma \cup \{ c_{k+1}<c_k<d_k<d_{k+1} :k<m\} \cup \{ \f _k^{\theta _k} :k\leq m\}$.
By (\ref{EQ205}) there is $n\in \o$ such that $[a_m,b_m]\subset J_n' \cong J_n \subset \X _{i_n}$.
Let $f: J_n' \rightarrow J_n$ be an isomorphism.
Then defining $c_k^\Y :=f(a_k)$ and $d_k^\Y :=f(b_k)$, for $k\leq m$,
and $\Y =\la X, \leq, c_0^\Y , \dots ,c_m^\Y,d_0^\Y , \dots ,d_m^\Y  \ra$
we have $c_m^\Y<\dots <c_0^\Y <d_0^\Y <\dots <d_m^\Y $.
For $k\leq m$ we have $[c_k^\Y ,d_k^\Y] \cong [a_k,b_k]$;
so, $[c_k^\Y ,d_k^\Y]\models \f _k$, that is $\Y \models \f _k^{\theta _k}$.

Let $\Z $ be a countable model of $\CT$.
Then the reduct $\Z ':=\Z \mid \{ \leq \}$ is a model of $\Th (\X)$
and by Fact \ref{T202}(c) $\Z '$ is a disjoint union of linear orders, $\Z'=\bcd _{s\in S}\Z_s$.
Since $\Z \models \CT$ we have $\dots <c_1^\Z <c_0^\Z <d_0^\Z <d_1^\Z<\dots $
and, hence, there is $s\in S$ such that $c_k^\Z ,d_k^\Z \in \Z _s$, for all $k\in \o$.
Suppose that $\s :=\tp (\Z _s)\in \CC$; then $\rank (\s )=m$, for some $m\in \o$.
Since $\Z \models \CT$ we have $[c_k^\Z ,d_k^\Z] \cong [a_k,b_k]$
and, hence $\bigcup _{k\in \o}[c_k^\Z ,d_k^\Z]$ is an interval in $\Z _s$ isomorphic to $J'$
and, since $\tp (J_n')=\tp (J_n)=\s _n$,
$\Z_s$ contains an interval of type $\s _n$.
So, by (i) and Fact \ref{T206}(b), $n=\rank (\s _n)\leq 2m+1$, for all $n\in \o$,
which is false.
Thus the component $\Z_s$ of $\Z '$ is not $\o$-categorical,
by Rubin's theorem we have $I(\Z _s)=\c$
and, since $\Z _s$ is of finite diameter, by  Proposition \ref{T201} we have $I(\X )=I(\Z ')=\c$.

{\it Subcase 2.2}: $\neg$ Subcase 2.1.
First let $n_0\in \o$, where $J_{n_0}'$ is an initial part of $J'$.
If there is $c=\min J_{n_0}'$,
then we put $a_n =c$, for all $n\in \o$ and work as above.
If $\min J_{n_0}'$ does not exist,
then we take a coinitial sequence $\la a_n :n\in \o\ra$ as above,
and the construction works,
because again $\Z_s$ contains intervals of type $\s _n$, $n\in \o$.
If there is $n_0\in \o$ such that $J_{n_0}'$ is a final part of $J'$,
we have a dual proof.
\kdok
By Theorem \ref{T203} if $\X =\bcd _{i\in \o} \X _i$, where $\X _i$ are $\o$-categorical linear orders,
then $I (\X )=1$ if, for example, $\X _i \cong \Q$, for all $i\in \o$, (then $\X$ is one of the ultrahomogeneous posets from the Schmerl list)
or if $\X _i$'s are finite linear orders of size $\leq n$, for some $n\in \N$.
On the other hand $I (\X )=\c$, if $\X =\bcd _{n\in \N}n$; then $\X$ is the atomic model of $\Th (\X )$,
the connectivity components are discrete linear orders with end points,
and the saturated model of $\Th (\X )$ is obtained by adding the union of $\o$-many copies of the linear order $\o + \zeta \eta +\o ^*$ to $\X$
(where $\zeta :=\tp (\Z)$ and $\eta :=\tp (\Q)$).
Thus the theory $\CT:=\Th (\bcd _{n\in \N}n )$ is small ($\bigcup _{n\in \N}S_n (\CT )\leq \o$).
\begin{ex}\label{EX000}\rm
{\it A large theory of a countable union of $\o$-categorical linear orders.}
Let $\Q=\{ q_n:n\in\N\}$ be an enumeration of the set of rational numbers
and $I:= \{ \la F,G\ra \in [\N]^{<\o}\times [\N]^{<\o}: F\cap G=\emptyset\}$.
For $\la F,G\ra \in I$ the linear order $\X _{F,G}:=\s (F) +1 +\s (G)$ is $\o$-categorical,
and for $\X :=\bcd _{\la F,G\ra \in I}\X _{F,G}$ we show that the theory $\CT :=\Th (\X )$ is large.

For $n\in \N$, an {\it $n$-jump} in a linear order $\Y$
is a sequence of its consecutive elements $y_0 <\dots <y_n$
such that $y_0$ does not have an immediate predecessor and  $y_n$ does not have an immediate successor.
Let $\p _n^- (x)$ be a formula saying that there is an $n$-jump below $x$; namely, defining $\f (u,v):= u<v \land \neg \exists w \;(u<w<v)$, let
$$\textstyle
\p _n^- (x):= \exists v_0,\dots,v_n \Big( \bigwedge _{i<n}\f (v_i,v_{i+1})\land v_n <x \land \neg\exists u \, \f (u,v_0)\land \neg\exists u \, \f (v_n,u)\Big)
$$
and, dually, let $\p _n^+ (x)$ be a formula saying that there is an $n$-jump above $x$.
For an irrational number $a$ let $M_a^-:=\{ n\in \N : q_n <a\}$ and $M_a^+:=\{ n\in \N : q_n >a\}$; clearly, $\{ M_a^- ,M_a^+\}$ is a partition of $\N$.
We show that the partial type
$$
\Pi _a (x):=\{ \p _n^- (x)\land \neg \p _n^+ (x) : n\in M_a^- \}\cup \{ \p _n^+ (x)\land \neg \p _n^- (x) : n\in M_a^+ \}
$$
is finitely satisfiable in $\X$.
So, if $F\in [M_a^- ]^{<\o}$ and $G\in [M_a^+ ]^{<\o}$,
then $\la F,G\ra \in I$
and the finite type  $\{ \p _n^- (x)\land \neg \p _n^+ (x) : n\in F \}\cup \{ \p _n^+ (x)\land \neg \p _n^- (x) : n\in G \}$
is satisfied in $\X$ by the ``middle" element of the component $\X _{F,G}:=\s (F) +1 +\s (G)$ of $\X$.
Thus there is $p_a (x)\in S_1 (\CT)$ such that $\Pi _a (x)\subset p_a (x)$.

If $a,b \in \BR \setminus \Q$ and $a<b$,
then there is $n\in \N$ such that $a <q_n <b$;
thus $n\in M_a^+  \cap M_b^-$
and, hence, $\p _n^+ (x)\in \Pi _a (x)$ and  $\neg \p _n^+ (x)\in \Pi _b (x)$,
which gives $p_a (x)\neq p_b (x)$.
So, $p_a (x)$, $a\in \BR \setminus \Q$, are different 1-types of $\CT$ and the theory $\CT$ is large indeed.
\end{ex}

\paragraph{Acknowledgement.}
This research was supported by the Science Fund of the Republic of Serbia,
Program IDEAS, Grant No.\ 7750027:
{\it Set-theoretic, model-theoretic and Ramsey-theoretic
phenomena in mathematical structures: similarity and diversity}--SMART.

\end{document}